\documentclass{article}
\usepackage[a4paper]{geometry}
\usepackage{amsmath}
\usepackage{amssymb}
\usepackage{graphicx}
\usepackage{subfig}

\usepackage{amsthm}
\usepackage{amsfonts}
\usepackage{amssymb}
\usepackage{mathrsfs}
\usepackage{mathtools}
\usepackage{xcolor,soul}
\usepackage{tikz}
\usepackage{pifont}

\usepackage{mathtools}

\usetikzlibrary{decorations.pathreplacing}

\newcommand{\eps} {\varepsilon}
\newcommand{\e} {\mathrm{e}}
\newcommand{\ii} {\mathrm{i}}
\newcommand{\diff}[2]{\frac{\mathrm{d} #1}{\mathrm{d} #2}}
\newcommand{\sign}{\mathrm{sign}}

\title{Capturing the cascade: a transseries approach to delayed bifurcations}
\author{In\^{e}s Aniceto$^1$\footnote{Corresponding Author. Electronic address: i.aniceto@soton.ac.uk} , Daniel Hasenbichler$^1$, Christopher J. Howls$^1$, Christopher J. Lustri$^2$}
\date{%
    $^1$Mathematical Sciences, University of Southampton, Highfield, Southampton, SO17 1BJ, UK\\[2ex]%
    $^2$Department of Mathematics and Statistics, 12 Wally's Walk, Macquarie University, New South Wales 2109, Australia
}         

\begin{document}

\maketitle

\abstract{Transseries expansions build upon ordinary power series methods by including additional basis elements such as exponentials and logarithms. Alternative summation methods can then be used to ``resum'' series to obtain more efficient approximations, and have been successfully widely applied in the study of continuous linear and nonlinear, single and multidimensional problems.
 In particular, a method known as transasymptotic resummation can be used to describe continuous behaviour occurring on multiple scales without the need for asymptotic matching.
 Here we apply transasymptotic resummation to discrete systems and show that it may 
 be used to naturally and efficiently describe discrete delayed bifurcations, or ``canards'', in singularly-perturbed variants of the logistic map which contain delayed period-doubling bifurcations. We use transasymptotic resummation to approximate the solutions, and describe the behaviour of the solution across the bifurcations. This approach has two significant advantages: it may be applied in systematic fashion even across multiple bifurcations, and the exponential multipliers encode information about the bifurcations that are used to explain effects seen in the solution behaviour.}

\section{Introduction}

Transseries are a natural extension to classical asymptotic power series which are used to study systems in which the solution behaviour depends on multiple distinct exponential scales. A transseries represents the solution to a system as the sum of multiple power series, each multiplied by a different exponential prefactor \cite{Edgar2009}. The value in this approach is that a transseries developed in one region of parameter space can typically be extended into regimes in which the solution depends on different scales, simply be applying different summation methods to the transseries itself, without the need to rebalance the equation terms and apply matched asymptotic expansions to connect the regimes.

Transseries and summation techniques have been used to study the behaviour of a wide range of parameter-dependent continuous systems. Applications include the study of general nonlinear ordinary differential equations  \cite{Costin1995,Costin1998,Howls2010,OldeDaalhuis2005a}, the first Painlev\'{e} equation \cite{Aniceto:2011nu,Garoufalidis2011,OldeDaalhuis2005}, topological string theory \cite{Couso-Santamaria:2014iia,Grassi:2014cla}, field theory and semi-classical quantum mechanics \cite{Aniceto:2013fka,Basar:2015xna,Demulder:2016mja,Dunne:2014bca,Grassi:2018spf,Marino:2020dgc}, relativistic hydrodynamics and Einstein partial differential equations \cite{Aniceto:2018uik,Casalderrey2018,Heller:2015dha}, and $q-$series and knot invariants \cite{Dorigoni:2020oon,Garoufalidis:2020nut}. More recently transseries methods have been been extended to study of discrete problems, such as particular matrix models governed by the first discrete Painlev\'{e} equation \cite{Aniceto:2011nu,Couso-Santamaria:2015wga,Marino:2008vx,Pasquetti:2009jg,Schiappa:2013opa}. The role of exponential scales in discrete maps and chaos has been previously analysed in the context of Stokes phenomena \cite{Gelfreich2001,Hakim1993,Martin2011,Shudo2008}.

The transasymptotic method introduced in \cite{Costin1998,Costin2001,Costin2015} consists of constructing a transseries in terms of some small parameter $\varepsilon$. The transseries terms are then reordered, and higher-order exponential terms at each order of the small parameter are summed. This change of summation order, or ``resummation'', captures the behaviour of the system in regions where different exponential terms dominate the solution. Transasymptotic methods have been been used to determine the location of moveable poles in Painlev\'{e} equations directly from asymptotic solutions \cite{Aniceto:toapp,Costin2015} (see also \cite{Costin:2019xql,Costin:2020pcj}). In this work, we show that these transseries and summation approaches can be used in a systematic fashion to identify complicated bifurcations in discrete systems that typically require careful application of multiple scales \cite{Hall2016} or renormalisation methods \cite{Baesens1991a}. 

 In this study, we demonstrate that transseries resummation can be used systematically to accurately capture the behaviour of the solution to a discrete equation containing a periodic-doubling cascade with delayed bifurcations.  Delayed bifurcations may occur in dynamical systems where an underlying parameter is itself slowly varying and the solution initially clings to a metastable branch of the solution before eventually jumping to the stable branch.  They have been studied widely in systems of ordinary differential equations (see, for example, \cite{Wechselberger2012}). These ``slow-fast'' systems have behaviour occurring on two (or more) distinct timescales, with the solution trajectory remaining near to an unstable solution for a significant distance after stability has been lost; solutions containing this behaviour have been termed ``canards''. 

There has been a significant volume of work studying the asymptotic behaviour of canards in continuous settings, see for example, using composite asymptotic expansions in \cite{Benoit1998,Fruchard2003,Eckhaus1983}, and Borel summation methods in \cite{dEscurac2017}. Borel summation methods are closely connected to transseries resummation methods (see \cite{Aniceto1999}), and have been used to study discrete problems, as in \cite{OldeDaalhuis2004}. This motivates the idea that transseries resummation methods could be a useful technique for studying delayed bifurcation behaviour; in this study, we will focus on delayed bifurcations appearing in discrete systems, and in particular, singularly perturbed variants of the logistic map. 

We will show that period doubling bifurcations depend on the interaction between different exponential factors, and it is therefore advantageous to represent them explicitly using transseries. By expanding in the asymptotic limit, we may determine terms in the algebraic power series to determine the initially stable non-periodic solution. The next step will be to reorder the transseries terms and perform a transasymptotic resummation, which will produce an accurate description of the doubling phenomena. This approach has the additional advantage that it allows us to determine further subdominant exponential scales in the transseries explicitly which dictate subsequent doubling bifurcations present in the solution.

By incorporating a multiple scales ansatz into the transseries expression, we will shown that transseries resummation -- which was developed to describe continuous behaviour -- can be used to calculate discrete variation without any further analysis to the transseries method. 

We study here two variants of the ubiquitous (and generic) standard logistic map 
\begin{equation}
y(n+1)=\lambda y(n)\left[1-y(n)\right],\qquad0<y(0)<1,\label{eq:Logistic-map0}
\end{equation}
where $\lambda$ is a dimensionless bifurcation parameter $0<\lambda\le4$.  

This system contains a period-doubling route to chaos, found by allowing the parameter $\lambda$ to vary. In the range $1 < \lambda <  3$, this system tends to a stable equilibrium without periodic effects.  In the range $3 < \lambda < 1 + \sqrt{6}$, the system tends to a 2-periodic stable equilibrium. Increasing $\lambda$ beyond $1 + \sqrt{6}$ produces systems that tend to stable equilibria with higher periodicity. 

The earliest study of the delayed bifurcations in the slowly-varying logistic map is \cite{Baesens1991a}, who applied renormalisation methods to derive asymptotic scaling laws for the delays between period doubling, and performed analysis and numerical experiments to determine the location of the bifurcation points. In addition to establishing specific results about the slowly-varying logistic map, this study established that delayed bifurcations can play an essential role in the behaviour of discrete systems. Similar methods were used to study delayed bifurcations in more general unimodal maps \cite{Davies1997}, as well as discrete maps with noise \cite{Baesens1991b, Davies2001, Davies2002}.

Further studies of this system appeared in subsequent years. In \cite{ElRabih2003,Fruchard1991,Fruchard2003}, the existence of canard solutions was rigorously proven in general classes of discrete maps that include the slowly-varying logistic map. Further discussions of canard solutions to both discrete continuous and discrete dynamical systems are given in \cite{Fruchard1996, Fruchard2009}.

In more recent years, this system was studied using matched asymptotic expansions and multiple scales methods \cite{Hall2016}. The purpose of this previous work was to show that the method of multiple scales could be used to combine a ``fast'' discrete timescale with a slow time variable that could be treated as continuous, while still capturing the essentially discrete-scale behaviour present in the problem. By carefully balancing terms, the authors were able to identify the bifurcation points and produce accurate asymptotic approximations to the solution behaviour on both sides of the delayed bifurcation. In the works described here, the slowly-varying logistic equation has provided a useful testing ground for treatments of discrete systems, due to the complicated behaviour that it produces.  

It has been shown in \cite{Howls2010} that transseries approaches may be used to improve upon asymptotic results obtained using matched 
asymptotic expansions. In that study, transseries resummation methods were used to obtain a uniform approximation to a continuous problem 
that had been previously solved using multiple scales methods. The transseries approach was able to naturally incorporate higher 
exponential terms, and thereby improve on the accuracy of the results, even for values of the perturbation parameter that were not 
extremely small. Motivated by this result, we will show that transseries resummation can be used to improve on existing multiple scales 
results in discrete systems. 

In section 2 we will first study the ``static'' logistic map, corresponding to $\lambda = 3 + \eps$ for a fixed choice of (small) real parameter $\eps$. This will be used to introduce and to demonstrate the technical details of the transseries process

\begin{figure}[tb]
\centering

\includegraphics[width=0.6\linewidth]{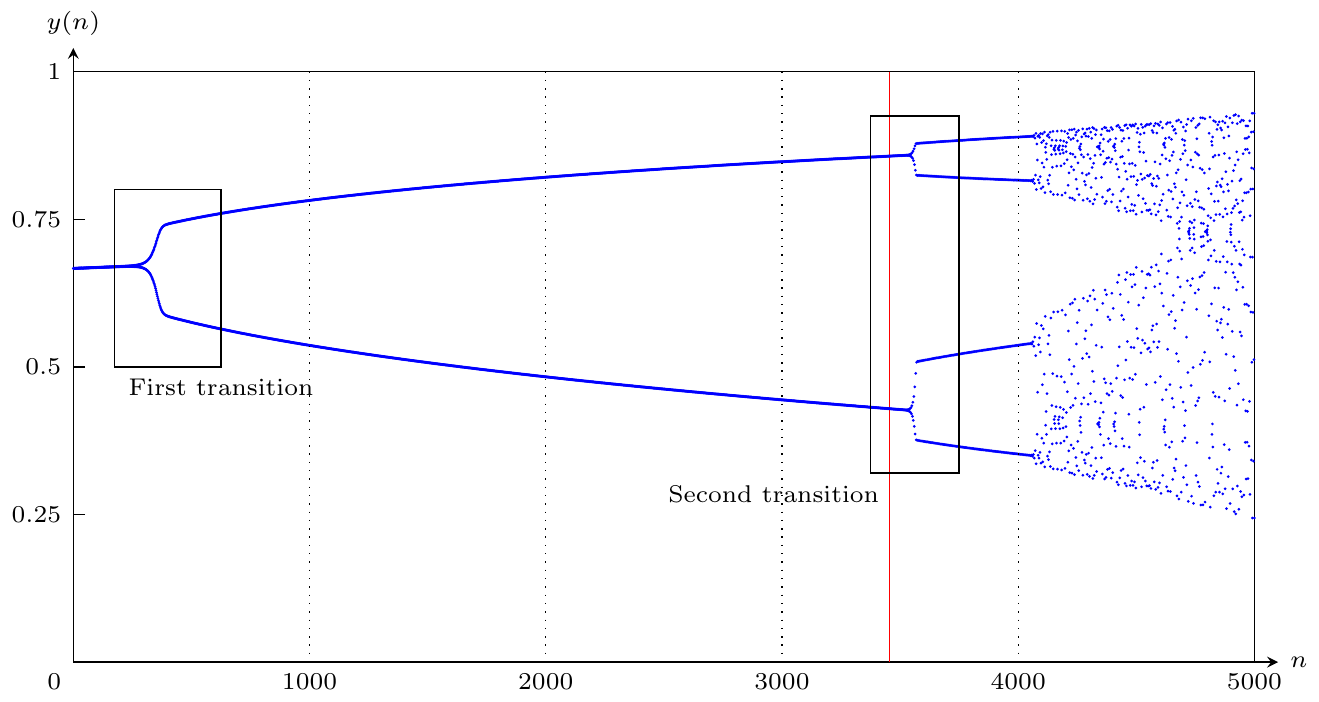}
\caption{Solution to logistic equation \ref{eq:Logistic-map0}, where $\lambda = 3 + \eps n$ with $\eps = 0.012^2$. The period-doubling cascade is apparent; the transition between non-periodic and 2-periodic behaviour is visible, as is the transition between 2-periodic and 4-periodic behaviour. As the solution continues, it eventually becomes chaotic. The 2-periodic behaviour in the solution begins to contribute immediately, but is not immediately visibly apparent due to the delay in the bifurcation behaviour. Similarly, from the analysis in Section \ref{s:dynamic}, we will determine that solution begins to display 4-periodic behaviour at $n \approx 3455$, shown as a red line, but it is also not immediately visibly apparent. }\label{fig:dyn-exact}
\end{figure}

The second variant, studied in section 3, is the ``dynamical'' logistic equation, which has a slowly increasing bifurcation parameter $\lambda = 3 + \eps n$. In this case, the bifurcation parameter grows, with different solutions becoming stable and unstable as $n$ increases. In Figure \ref{fig:dyn-exact},

Apart from the pedagogical aspects of demonstrating the applicability of transseries to such problems to recover and extend existing results, in this work we demonstrate four main enhancements.  

Firstly, we show that transseries resummation can capture the behaviour of solutions to both a static and 
slowly-varying logistic map in a systematic and efficient fashion. The multiple scales method used by \cite{Hall2016} was able 
to produce asymptotically valid approximations to the solution, but the process requires a careful 
expansion and asymptotic matching each time a bifurcation occurs. The transseries resummation approach used here 
is systematic, without any need for matching, and can be applied in largely identical fashion to capture each successive bifurcation. 

Secondly, transseries resummation allows us to capture behaviour when the bifurcation parameter is not necessarily small. In 
\cite{Hall2016}, the static logistic map was studied in the limit that the bifurcation parameter was close to three, leading to 2-periodic
behaviour in the static map. In this study, the transseries resummation approach allows for the study of larger values of the bifurcation parameter, 
describing 4- and 8-periodic behaviour in the static map.

Thirdly, we demonstrate that transseries can capture and control the onset of period-doubling behaviour through in terms of exponential 
weights in the transseries coupled with their resummation.  By studying these exponential terms, we are able to determine precisely when 
higher-periodicity behaviour appears in the solution, and when the bifurcation starts to grow.  We will show that transseries resummation can approximate the behaviour of 2- and 4-periodic solutions, and that calculating the transseries exponential terms can explain the onset of higher  4- or even 8-periodic behaviour as the bifurcation parameter grows.  In principle the method could be continued in the same fashion to determine this behaviour as the dynamic map sweeps through subsequent bifurcations.

Fourthly, we are able to use transseries summation to significantly improve the approximation accuracy of the solutions over and above that afforded by matched asymptotics in several parameter regimes, with minimal, if any, additional effort. The increase in accuracy is a consequence of the inclusion of multiple exponential scales in the solution approximation, and is most apparent in parameter regimes in which different exponential scales all contribute to the solution behaviour.

\section{Static Logistic Equation}

First we consider the static logistic equation, given by
\begin{equation}
y(n+1)=(3 + \eps) y(n)\left[1-y(n)\right],\qquad y(0) = 2/3.\label{eq:Logistic-map}
\end{equation}
We will write the solution as a continuous transseries in terms of $\eps>0$. We first produce an asymptotic expansion in the limit $0<\eps\ll1$ but as we shall see later, the transseries approach will be used to extend this result to produce an accurate approximation for $\eps=\mathcal{O}(1)$. We will then show that this continuous transseries is capable of capturing discrete period-doubling effects seen in this system, and approximating higher periodicity behaviour for values of $\eps$ that lead to 2-, 4- and even 8-periodic solutions.

In \cite{Hall2016}, the authors studied the asymptotic behaviour of this system for small $\varepsilon$ using multiple scales methods. This showed the manner in which the behaviour approached the 2-periodic stable manifold associated with $\lambda > 3$. Using transseries methods, we can extend this approach to consider systems in which $\varepsilon$ is not asymptotically small, and demonstrate the manner in which the solution approaches the stable solution for higher periodicities. 

In order to first determine the non-periodic and 2-periodic solutions, we ignore the initial condition and solve \eqref{eq:Logistic-map} with the condition that $y(n+2) = y(n)$. This gives three unique solutions. One solution is non-periodic, and is given by
\begin{equation}
y(n) = \frac{2+\eps}{3+\eps}.
\end{equation}
The remaining two solutions are 2-periodic, and are given by
\begin{equation}
y(n) = \frac{4 + \eps \pm (-1)^n\sqrt{\eps(4+\eps)}}{2(3+\eps)}.\label{eq:2per-manifold}
\end{equation}
For $\eps > 0$, the non-periodic solution is unstable. For $0 < \eps < \sqrt{6}-2$, the 2-periodic solution is stable. If $\eps$ exceeds $\sqrt{6}-2$, the 2-periodic solution is unstable, and the stable solution to the system becomes 4-periodic, and can be identified by solving $y(n+4) = y(n)$, however this solution cannot be expressed in closed form. Continuing to increase $\eps$ leads to the periodicity of the stable solution increasing until chaotic behaviour is eventually obtained.

\subsection{2-periodic solution}\label{S:static-2-per}

\subsubsection{Transseries ansatz}

We begin by applying a transseries ansatz, including a continuous variable $x$ and the small parameter $\varepsilon$. We set $x = \varepsilon n$, assuming for now $\eps\ll 1$, and $R(x) = y(n)$ obeys
\begin{equation}
R(x+\varepsilon)=(3+\varepsilon)R(x)[1-R(x)],\qquad R(0)=\tfrac{2}{3}.\label{eq:logistic-finite-diff-constant}
\end{equation}
We formulate an ansatz for the solution behaviour in terms of both $\varepsilon$ and a transseries parameter $\sigma_0$. The preliminary ansatz is given as the standard (non-logarithmic) transseries:
\begin{equation}
R(x,\eps; \sigma_0) = \sum_{m=0}^{\infty} \sigma_0^m \e^{-mA(x)/\eps} R_m(x,\eps),\qquad R_m(x,\varepsilon) = \eps^{\beta_m}\sum_{k=0}^{\infty}\varepsilon^k R_{m,k}(x),\label{eq:logistic-transseries-constant}
\end{equation}
where $\beta_m$ will be chosen such that $R_{m,0}$ takes a nonzero value. Applying this expression to \eqref{eq:logistic-finite-diff-constant} allows for the system to be matched in powers of $\sigma_0$. At leading order, we find
\begin{equation}
R_0(x,\eps) = \frac{2+\varepsilon}{3+\varepsilon}.\label{eq:logistic-R0}
\end{equation}
This expression gives the non-periodic manifold, which is stable for $1 < \lambda < 3$. For values of $\lambda$ greater than three, corresponding to positive $\varepsilon$, we expect to see additional periodic effects emerge from the transseries expression. Continuing to the next order in $\sigma_0$, we find that
\begin{equation}
\e^{-[A(x+\eps) - A(x)]/\eps} R_1(x + \eps,\eps) = -(1+\eps) R_1(x,\eps).
\end{equation}
By expanding $R_1$ as a power series in $\eps$, \eqref{eq:logistic-transseries-constant} gives
\begin{equation}
\e^{-[A(x+\eps) - A(x)]/\eps}  \sum_{k=0}^{\infty}\varepsilon^k R_{1,k}(x) = -(1+\eps)  \sum_{k=0}^{\infty}\varepsilon^k R_{1,k}(x).\label{eq:logistic-R1}
\end{equation}
Expanding as a Taylor series in $\varepsilon$ gives $A(x + \varepsilon) = A(x) + \varepsilon A'(x) + \cdots$. 
 Matching equal powers of $\varepsilon$ on both sides in \eqref{eq:logistic-R1} therefore requires $\mathrm{exp}(-A'(x)) = -1$, or
\begin{equation}
A(x) = (2p+1)\pi\ii  x,\qquad p \in \mathbb{Z}.\label{eq:logistic-A}
\end{equation}
The arbitrary constant in $A(x)$ may be absorbed into the transseries parameter $\sigma_0$, and is therefore set to be zero for convenience. The solution will only be evaluated for $n \in \mathbb{Z}$, which corresponds to $x/\varepsilon \in \mathbb{Z}$. The choice of $p$ has no effect at these values, and we therefore set $p = 0$ without loss of generality.

Using the expressions for $R_0$ in \eqref{eq:logistic-R0} and $A$ in \eqref{eq:logistic-A} in \eqref{eq:logistic-finite-diff-constant} and tracking orders of $\sigma_0$, we obtain a recurrence relation for $R_m$:
\begin{align}
(-1)^{m}R_{m}(x+\varepsilon,\eps)  =-(1+\eps)R_{m}(x,\eps)-(3+\eps)\sum_{j=1}^{m-1}R_{j}(x,\eps)R_{m-j}(x,\eps),\label{eq:recursion}
\end{align}
for $m \geq 1$, where we take the convention that the summation term is zero when $m = 1$. It is straightforward to show by direct substitution that the solution to this recurrence takes the form
\begin{equation}
R_m(x,\eps) = \e^{m x \log(1+\eps)/\eps} \overline{R}_m(\eps),\label{eq:Rm-general} \qquad \overline{R}_m(\eps) 
= \eps^{\beta_m} \sum_{k=0}^\infty
\eps^k \, \overline{R}_{m,k},
\end{equation}
where the functions $\overline{R}_m(\eps)$ can be computed directly. $\overline{R}_1(\eps)$ may be chosen to be an arbitrary constant, which again can be absorbed into $\sigma_0$. For algebraic convenience, we can select $R_1(\eps) = \eps$. The subsequent terms may hence be obtained using the recursion \eqref{eq:recursion}, which gives the first few terms as
\begin{equation}
\overline{R}_{2}(\eps)  =-\frac{(3+\eps)\eps^2}{(1+\eps)(2+\eps)},\qquad \overline{R}_{3}(\eps)  = \frac{2(3+\eps)^{2}\eps^{2}}{(1+\eps)^{2}(2+\eps)^{2}},\qquad \overline{R}_4(\eps) = -\frac{(3+\varepsilon)^{3}\left(\varepsilon-4\right)\varepsilon^{3}}{(1+\varepsilon)^{3}(2+\varepsilon)^{3}\left(1+\varepsilon+\varepsilon^{2}\right)}.\label{eq:logistic-early-terms}
\end{equation}
By analysing the form of the recurrence solution, it is straightforward to determine $\beta_m$ in general, giving $\beta_m = (m+1)/2$ for $m$ even, and $\beta_m = (m+2)/2$ for $m$ even. It can also be seen by direct calculation that
\begin{equation}
\overline{R}_{2m+1,0}=\frac{(-9)^{m}\Gamma\left(m+\frac{1}{2}\right)}{\sqrt{\pi}\Gamma\left(n+1\right)}, \qquad \overline{R}_{2m,0}=\frac{(-9)^{m}}{6}.\label{eq:logistic-eq-exp-coeffs-constant}
\end{equation}

\subsubsection{Computing the terms in the resummed transseries}

The alternating behaviour of \eqref{eq:logistic-eq-exp-coeffs-constant} combined with the general form of $R_m(x, \eps)$ described in \eqref{eq:Rm-general} suggests that the ansatz may be conveniently re-written to incorporate these elements explicitly. In particular, the fact that we have two different sub-series in (\ref{eq:logistic-eq-exp-coeffs-constant}) depending
on the parity of $m$ 
suggests that we should split the series into odd and even powers of $m$. We therefore write the ansatz \eqref{eq:logistic-transseries-constant} as
\begin{align}
R(x,\eps;\sigma_0) = R_0(\eps) + \sqrt{\eps}\sum_{k=0}^{\infty} \eps^k \sum_{m=0}^{\infty} \left(\sigma_0 \sqrt{\varepsilon} \mathrm{e}^{-A(x)/\varepsilon+x\log(1+\varepsilon)/\varepsilon}\right)^{2m+1}
\overline{R}_{2m+1,k}&\nonumber \\
 +\varepsilon\sum_{k=0}^{\infty}\varepsilon^{k}\sum_{m=1}^{\infty}\left(\sigma_0 \sqrt{\varepsilon} \mathrm{e}^{-A(x)/\varepsilon+x\log(1+\varepsilon)/\varepsilon}\right)^{2m}\overline{R}_{2m,k}&,\label{eq:2per-general-ansatz}
\end{align}
where we have switched the order of summation, noting that the exponential terms are all of the same order in $\eps$. Note that in the static case both sums in \eqref{eq:2per-general-ansatz} are convergent. We define a new series parameter $\tau_0$, as well as odd and even power series in this parameter, such that
\begin{equation}
\tau_0(x,\varepsilon)=\sigma_0\sqrt{\varepsilon}\mathrm{e}^{-A(x)/\varepsilon+x\,\log(1+\varepsilon)/\varepsilon},\qquad \Omega_{o,k}(\tau_0) = \sum_{m=0}^{\infty}\tau_0^{2m+1}\overline{R}_{2m+1,k},\qquad \Omega_{e,k}(\tau_0) = \sum_{m=0}^{\infty}\tau_0^{2m}\overline{R}_{2m,k}. \label{eq:logistic-tau-constant}
\end{equation}
The transseries expression is now given by (the  $x,\sigma_0$ dependence is encoded in $\tau_0$):
\begin{align}
R(\tau_0,\eps) = R_0(\eps) + \sqrt{\eps} \sum_{k=0}^{\infty}\eps^{k} \Omega_{o,k}(\tau_0) + \eps \sum_{k=0}^{\infty} \eps^k \Omega_{e,k}(\tau_0).\label{eq:logistic-R-newdef}
\end{align}
 We may apply this "resummed" transseries to the logistic equation \eqref{eq:logistic-finite-diff-constant} and equate powers of $\eps$ to obtain expressions for $\Omega_{o,k}$ and $\Omega_{e,k}$. This process is somewhat technical, and is shown explicitly in Appendix \ref{A:logistic-transseries}. Equating terms of order $\eps$ and $\eps^{3/2}$ respectively gives
\begin{equation}
\Omega_{e,0} = -\frac{3}{2}(\Omega_{o,0})^2,\qquad \tau_0 \diff{\Omega_{o,0}}{\tau_0} = \Omega_{o,0} - 9(\Omega_{o,0})^3.\label{eq:OmegaSys}
\end{equation} 
Solving the ordinary differential equation gives $\Omega_{o,0}(\tau_0) = \pm \tau_0 ( C + 9\tau_0^2)^{-1/2}$. The sign and constant may be chosen arbitrarily, as this choice again can be absorbed into $\sigma_0$, which is yet to be determined. In order to maintain consistency with \eqref{eq:logistic-early-terms}, we select the positive sign and $C = 1$, giving
\begin{align}
\Omega_{o,0}(\tau_0)  =\frac{\tau_0}{\sqrt{1+9\tau_0^{2}}},\qquad \Omega_{e,0}(\tau_0)  =-\frac{3}{2}\frac{\tau_0^{2}}{(1+9\tau_0^{2})}.\label{eq:static-omega0}
\end{align}
Continuing this process and equating terms of order $\eps^{2}$ and $\eps^{5/2}$ respectively gives
\begin{align}
\Omega_{e,1}  =\frac{\tau_0^{2}\left(14-9\tau_0^{2}\right)}{8(1+9\tau_0^{2})^{2}}-3\Omega_{o,0}\Omega_{o,1},\qquad 
\tau_0\diff{\Omega_{o,1}}{\tau_0}  =\frac{1-18\tau_0^{2}}{1+9\tau_0^{2}}\Omega_{o,1}+\frac{3\tau_0^{2}(28-45\tau_0^{2})}{4\left(1+9\tau_0^{2}\right)^{5/2}}.
\end{align}
Solving these equations gives
\begin{align}
\Omega_{o,1}(\tau_0)  =-\frac{\tau_0(45\tau_0^{2}-33\log(1+9\tau_0^{2}))}{24(1+9\tau_0^{2})^{3/2}},\qquad \Omega_{e,1}(\tau_0) =\frac{\tau_0^{2}\left(14+36\tau_0^{2}-33\log\left(1+9\tau_0^{2}\right)\right)}{8(1+9\tau_0^{2})^{2}},\label{eq:static-omega1}
\end{align}
where yet again the arbitrary constant of integration can be chosen arbitrarily with the choice being absorbed into $\sigma_0$, and was therefore selected in order to maintain consistency with \eqref{eq:logistic-early-terms}. In principle, this process can be continued to higher transasymptotic order indefinitely by matching at higher orders of $\eps$.

We can compare these expressions to the multiple scales expansion obtained in \cite{Hall2016}. A straightforward comparison shows that the expression $\Omega_{o,1}$ gives the first correction to the composite expansion, denoted $P(s)$, from \cite{Hall2016}. The expression for $P(s)$ contains an exponential multiplier $\e^s$, which corresponds to the exponential scaling in $\tau_0$. This confirms that the resummed transseries identifies the region in which the power series expression breaks down, and allows for the expansion to be continued past this region even without the use of matched asymptotic expansions. 

\subsubsection{Initial value problem}

In order to complete the approximation, we must determine the value of $\sigma_0$ using the initial condition in \eqref{eq:logistic-finite-diff-constant}, which requires that $R(x=0,\eps;\sigma_0) = 2/3$. Noting that $\tau_0(x=0,\eps) = \sqrt{\eps}\sigma_0$, this corresponds to solving
\begin{equation}
R_0(\eps) + \sqrt{\varepsilon}(\Omega_{o,0}(\sigma_0\sqrt{\varepsilon})+\varepsilon\Omega_{o,1}(\sigma_0\sqrt{\varepsilon}))+\varepsilon(\Omega_{e,0}(\sigma_0\sqrt{\varepsilon})+\varepsilon\Omega_{e,1}(\sigma_0\sqrt{\varepsilon}))+\mathcal{O}(\varepsilon^{4}) = \tfrac{2}{3}.\label{eq:logistic-ic}
\end{equation}
By expanding $\sigma_0$ as a power series in $\eps$ such that
\begin{equation}
\sigma_0(\eps) = \sum_{j=0}^{\infty} \eps^j \sigma_{0,j}.
\end{equation}
Matching powers of $\eps$ in \eqref{eq:logistic-ic} now gives
\begin{equation}\label{eq:per2-sigma1}
\sigma_{0,0} = -\tfrac{1}{9},\qquad \sigma_{0,1} = \tfrac{4}{81},\qquad \sigma_{0,2} = -\tfrac{19}{648}.
\end{equation}
This process may be continued as $\Omega_{o,k}$ and $\Omega_{e,k}$ are computed for higher values of $k$, giving the increasingly complete expression for the 2-periodic behaviour. 

\subsection{4-periodic solution}\label{s:static-4per}

\subsubsection{Transseries ansatz}

To obtain the 4-periodic solution requires an adaptation of the previous process. We now take a four-periodic perturbation about the 2-periodic solution obtained in \eqref{eq:logistic-ic}. This allows us to form a transseries that can be used to capture solutions which tend to a four-periodic stable manifold. In \cite{Hall2016}, this would have required solving a challenging multiple scales problem, as the asymptotic solution obtained therein is only valid for small $\eps$. Using the transseries approach, we obtain a significant more general result.

We write the solution as a perturbation around the non-period behaviour and the 2-periodic behaviour captured by the transseries expression \eqref{eq:logistic-R-newdef}, in terms of the variable $\tau_0$, here written as $\hat{R}(\tau_0,\eps)$. 
\begin{equation}
R(x,\eps) = R_0(\eps) + \sqrt{\eps} \sum_{k=0}^{\infty}\eps^{k} \Omega_{o,k}(\tau) + \eps \sum_{k=0}^{\infty} \eps^k \Omega_{e,k}(\tau) + S(x,\eps) = \hat{R}(x,\eps) + S(x,\eps).
\end{equation}
We will then show that this perturbation starts contributing for values of $\eps$ large enough. We can see by direct substitution into \eqref{eq:logistic-finite-diff-constant} that
\begin{equation}
S(x + \eps,\eps) = (3+\eps)(1-2\hat{R}(\tau_0,\eps) -S(x,\eps))S(x,\eps).\label{eq:deltaR-eqn}
\end{equation}
In order to identify the correct scaling for $S(x)$, we note the form of the 2-periodic manifold, given in \eqref{eq:2per-manifold}. We represent the 2-periodic behaviour in terms of the continuous variable $x$ by writing $(-1)^n$ as $\alpha=\sign(\cos(\pi x/\eps))$, such that
\begin{equation}
\hat{R}(\tau_0,\eps) = \frac{4 + \eps - \alpha\sqrt{\eps(4+\eps)}}{2(3+\eps)} + \mathcal{O}(\tau_0^{-1}) \quad \mathrm{as} \quad x \rightarrow \infty, \label{eq:R-largex}
\end{equation}
where the asymptotic order of this expression can be obtained by rewriting \eqref{eq:logistic-ic} in powers of $\tau_0^{-1}$ and equating terms. We may now follow similar methods to the 2-periodic case, and formulate an ansatz for the solution in terms of $\eps$ and a new transseries parameter, denoted $\sigma_1$. Motivated by \eqref{eq:Rm-general}, we choose the ansatz
\begin{equation}
S(x,\eps) = \sum_{m=1}^{\infty}\sigma_1^m \e^{-m B(x,\varepsilon)/\eps} S_m(\eps), \label{eq:4per-Rsum}
\end{equation}
noting that we now allow an $\eps$ dependence on the exponential scale. The exponential scaling $B(x,\eps)$ may then be determined by considering the large-$x$ behaviour. This is convenient, as we have the form for $R(x)$ in this limit, given in \eqref{eq:R-largex}, and we know from the asymptotic order of this expression that the solution approaches this limit exponentially as $x$ becomes large. Using \eqref{eq:R-largex} in \eqref{eq:deltaR-eqn} and  
matching powers of $\varepsilon$ in an identical fashion to \eqref{eq:logistic-R1} gives 
\begin{equation}\label{eq:static-bd}
\frac{\partial}{\partial x}B(x,\eps) = -\pi \ii   - \log\left(1 -\alpha\sqrt{\eps(4+\eps)}\right).
\end{equation}
This may be solved to give
\begin{equation}
B(x,\eps) =f(\eps) x - \eps g(x,\eps),\label{eq:4per-A1}
\end{equation}
where 
\begin{equation}
    f(\eps) = -\tfrac{1}{2}\log(1-\eps(4+\eps)) - \pi \ii,
\end{equation}
and $g(x,\eps)$ is a bounded function that vanishes for $x = \eps n$ for $n \in \mathbb{Z}$, given by
\begin{equation}
g(x,\eps) = \alpha\left(\frac{x}{2\eps} - \frac{1}{2}\left\lfloor\frac{x}{\eps} + \frac{1}{2}\right\rfloor\right)\log\left(\frac{1 - \sqrt{\eps(4+\eps)}}{1 + \sqrt{\eps(4+\eps)}}\right).
\end{equation}
As $g(n\eps,\eps) = 0$ for $n \in \mathbb{Z}$, this expression could be ignored and the result will still give the correct value of $B(x,\eps)$, and hence the correct exponential scaling, on $x = \eps n$. This therefore suggests that we can capture the 4-periodic solution by defining a new variable $\tau_1$, such that
\begin{equation}
\tau_1(x,\eps) = \sigma_1\sqrt{\eps}\e^{-B(x,\eps)/\eps}.\label{eq:4per-tau1}
\end{equation}
In subsequent analysis, it will be useful to have a convenient expression for the value of $\tau_1$ at $x + \eps$ and $x + 2\eps$. Through direct substitution, we find that
\begin{align}
\tau_1(x+\eps,\eps) &= -\ii({\eps(4+\eps)-1})^{1/2}\e^{-2g(x,\eps)}\tau_1(x,\eps), \label{eq:4per-tau1-1jump} \\ 
\tau_1(x+2\eps,\eps) &=-(\eps(4+\eps)-1) \,\tau_1(x,\eps).\label{eq:4per-tau1-2jump}
\end{align}

\subsubsection{Exponential Weights}\label{S:where-4-periodic-behaviour-matters}

\begin{figure}
\centering

\includegraphics{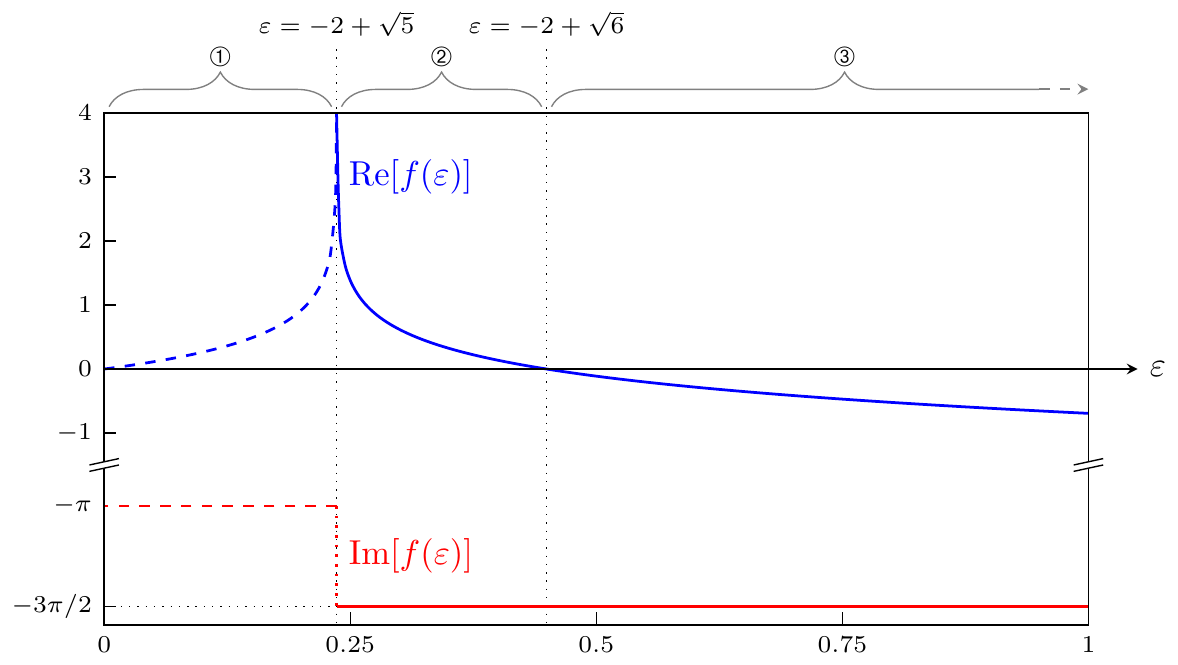}

\caption{This figure shows the real and imaginary parts of $f(\eps) = -\tfrac{1}{2}\log(1 - \eps (4+\eps)) - \pi \ii$, where the exponential weight $B(x)$ is given in \eqref{eq:4per-A1}. If $\mathrm{Re}[f(\eps)] > 0$, the 4-periodic exponential contribution is exponentially small in $\eps$, while if $\mathrm{Re}[f(\eps)] < 0$, the contribution is large, and must be incorporated into any approximation in order to accurately describe the system behaviour. In parameter regime \ding{192}, this exponential contribution is not present in the transseries, and is therefore denoted as a dashed curve. In regime \ding{193}, the 4-periodic exponential terms appear, but are exponentially small. In regime \ding{194}, the exponential contributions become large, and 4-periodic behaviour becomes apparent in the solution. The 4-periodicity of the solution arises due to $\mathrm{Im}[f(\eps)]$. This represents a multiplicative factor in $B(x)$ of $-\ii$, corresponding to 4-periodic behaviour in the exponential term.}\label{fig:A1-log}
\end{figure}

The form of $B(x,\eps)$ provides insight into the behaviour of the solution as $\eps$ increases outside of the range of validity of the original small-$\eps$ transseries. The behaviour of this term is shown in Figure \ref{fig:A1-log}, which identifies three distinct ranges of $\eps$ which must be considered separately. 

From the form of \eqref{eq:4per-tau1}, we see that $B(x,\eps)$ is the exponential controlling factor for $S$, and therefore determines how this series will contribute as $x$ grows. If $\mathrm{Re}[B] > 0$, corresponding to $\mathrm{Re}[f(\eps)] > 0$, the exponential contribution will decay as $x$ grows, while if $\mathrm{Re}[B] < 0$, corresponding to $\mathrm{Re}[f(\eps)] < 0$, the exponential part will grow and become the most significant contribution for large $x$.  

This change in sign occurs at $\eps = -2 + \sqrt{6}$. At this point, $\mathrm{Re}[f(\eps)]$ becomes negative, and the exponentials in \eqref{eq:4per-Rsum} therefore grow as $x$ becomes large, rather than decaying. This means that in Region \ding{194} the $S$ term is no longer a small decaying perturbation around $\hat{R}$, but rather plays a significant role in the limiting behaviour as $x \to \infty$. If $S$ is ignored, this behaviour is not captured in the transseries, and the resultant expression for $R$ is an inaccurate description of the solution behaviour. 

We note that $\mathrm{Im}[f(\eps)] = -3\pi/2$ for $\eps > -2 + \sqrt{5}$. This has the effect of making $\tau_1^m$ 4-periodic in $m \in \mathbb{Z}^+$, due to having a factor of $-\ii$, rather than the 2-periodic behaviour associated with a factor of $-1$. Hence, this exponential behaviour represented by $B(x,\eps)$ corresponds to 4-periodic effects in the solution. 

In order to include this behaviour in the transseries expression, we cannot simply expand the solution about $\eps = 0$. We must instead expand $S$ about some point $\eps_0$ such that the 4-periodic behaviour is present in the expansion.  This requirement suggests that $\eps_0 = -2 + \sqrt{6}$ is a sensible choice, as 4-periodic effects are apparent in the solution at this value.

Finally, we must consider the region in which the series terms obtained by expanding about $\eps_0$ are valid. If we examine the behaviour of $B(x,\eps)$ for $\eps > \eps_0$, we see that the real part of $f(\eps)$ becomes infinite as $\eps \rightarrow -2 + \sqrt{5}$. This corresponds to the exponentials disappearing, as every exponential term tends to zero. The series expansion about $\eps_0$ is not valid for $\eps$ smaller than this value. Consequently, region \ding{193} contains exponentially small 4-periodic behaviour, while no such behaviour exists in region \ding{192}. In Figure \ref{fig:A1-log} we have represented $f(\eps)$ in region \ding{192} as a dashed curve, to indicate that it does not have any effect on the transseries. 

Consequently, simply by studying $B(x,\eps)$, we are able to describe the onset of 4-periodicity in the solution. In region \ding{192}, there are no 4-periodic effects present. In region \ding{193}, there are 4-periodic effects caused by the appearance of new exponential terms, but they are exponentially subdominant compared to the 2-periodic behaviour. In region \ding{194}, these effects grow to become the most significant effect in the solution behaviour. Note that the switching of the 4-period exponentials is independent of the initial data, the latter only determines how quickly they grow to dominate the solution. For higher values of $\eps$, there must be values for which higher-periodicity behaviour appears. We discuss the onset of 8-periodic behaviour in Section \ref{S:8-per}.

Finally, we note that the change in exponential contribution has a parallel with a Borel transform approach to asymptotic expansions. Borel transforms encode the different exponential weights of an asymptotic series as singularities in a complex domain known as the ``Borel plane''. There, a change in the number of exponential contributions corresponds to singularities moving across a branch cut onto a different Riemann sheet of the Borel plane, giving rise to behaviour known as the ``Stokes phenomenon'' \cite{Berry1989} as the number of exponential contributions in an asymptotic series abruptly changes.  A similar, but not identical, behaviour occurs in this system at $\eps = -2 + \sqrt{5}$, where $\mathrm{Re}[f(\eps)]$ becomes infinite and $\mathrm{Im}[f(\eps)]$ changes instantaneously, corresponding to a branch point in the $f$-plane.  For more information on Borel transform methods, and their links to transseries and transasymptotic summations, see \cite{Aniceto1999,Costin1998,Costin2015,Howls2004,OldeDaalhuis2005a,OldeDaalhuis2005}.

\subsubsection{Computing the terms in the resummed transseries}

Writing an appropriate form for the 4-periodic ansatz is slightly more involved than in the 2-periodic case, given in \eqref{eq:2per-general-ansatz}. Recall from Section \ref{S:where-4-periodic-behaviour-matters} that the significant change in the behaviour of the exponential contribution occurs for values of $\eps$ greater than $\eps_0 = -2 + \sqrt{6}$. We therefore define a new series variable $\sqrt{6} \eta = \eps - \eps_0$, where the $\sqrt{6}$ term is included for subsequent algebraic convenience. 

In analogous fashion to \eqref{eq:logistic-R-newdef}, we again divide the ansatz up into separate power series. In the 2-periodic case, it was clear from the form of the previously calculated terms that splitting the odd and even powers of $\tau_0$ would capture the discrete variation effectively. From the analysis in Section \ref{S:where-4-periodic-behaviour-matters}, we determine that the power series for the 4-periodic solution should instead be split into four parts, such that
\begin{align}
S(\tau_1,\eta) = \sqrt{\eta}\sum_{k=0}^{\infty} \eta^k & \sum_{m=0}^{\infty}\tau_1^{4m+1} S_{4m+1,k}  + \eta \sum_{k=0}^{\infty}\eta^{k}\sum_{m=0}^{\infty}\tau_1^{4m+2}S_{4m+2,k}\nonumber \\
&+\sqrt{\eta}\sum_{k=0}^{\infty} \eta^k \sum_{m=0}^{\infty}\tau_1^{4m+3} S_{4m+3,k}  + \eta \sum_{k=0}^{\infty}\eta^{k}\sum_{m=0}^{\infty}\tau_1^{4m+4}S_{4m+4,k}.\label{eq:4per-general-ansatz}
\end{align}
Consequently, we now write each split power series as functions $\Theta_{j,k},\,j=1,2,3,4$, giving 
\begin{align}
S(\tau_1,\eta) = \sqrt{\eta}\sum_{k=0}^{\infty} \eta^k & \Theta_{1,k}(\tau_1)  + \eta \sum_{k=0}^{\infty}\eta^{k}\Theta_{2,k}(\tau_1) +\sqrt{\eta}\sum_{k=0}^{\infty} \eta^k \Theta_{3,k} (\tau_1)  + \eta \sum_{k=0}^{\infty}\eta^{k}\Theta_{4,k}(\tau_1) .\label{eq:4per-general-ansatz-2}
\end{align}
Noting that each series consists only of powers $\tau_1^m$ with the same $m\mod 4$, and comparing this with the expression for $\tau_1$ in \eqref{eq:4per-tau1} indicates that the functions $\Theta_{j,k}$ for $j = 1,\ldots, 4$ must have the symmetries
\begin{align}
\Theta_{1,k}(-\ii\tau_1) = -\ii\Theta_{1,k}(\tau_1),& \qquad \Theta_{2,k}(-\ii\tau_1) =  -\Theta_{2,k}(\tau_1), \label{eq:4per-Theta12}\\
\Theta_{3,k}(-\ii\tau_1) =  \ii\Theta_{3,k}(\tau_1), &\qquad \Theta_{4,k}(-\ii\tau_1) =  \Theta_{4,k}(\tau_1). \label{eq:4per-Theta34}
\end{align}

At this stage, it might be expected that we should express the governing equation \eqref{eq:deltaR-eqn} in terms of $\eta$, and perform an expansion in this variable; however, a comparison of the terms in \eqref{eq:4per-tau1-1jump} and \eqref{eq:4per-tau1-2jump} suggests that iterating the map once leads to a simplification. Writing the $x$ dependence explicitly, the equation becomes:
\begin{align}
    S(x+2\eps,\eps)= (3+\varepsilon)^2 [1-2\hat{R}(x+\varepsilon,\eps)-(3+\varepsilon)(1-& 2\hat{R}(x,\eps)-S(x,\eps))S(x,\eps)]\nonumber\\
   &\times(1-2\hat{R}(x,\eps)-S(x,\eps))S(x,\eps),\label{eq:4per-double-gov-eqn}
\end{align}
This expression does not contain any $S(x + \eps,\eps)$ terms, and instead only contains the double iteration term, $S(x + 2\eps,\eps)$. This is convenient, as the expression for $\tau_1(x+2\eps,\eps)$ is substantially simpler than $\tau_1(x+\eps,\eps)$, as it does not contain $g(x,\eps)$. This simplifies significantly the subsequent analysis. 

Expressing the left-hand term in \eqref{eq:4per-double-gov-eqn} in terms of $\tau_1$ and $\eta$ gives
\begin{align}
S(x+2\eps,\eps)= S(-(\eps(4+\eps)-1) \,\tau_1,\eps) = S(-(1 + 12\eta + 6\eta^2)\tau_1,\eta).
\end{align}
Rewriting \eqref{eq:4per-double-gov-eqn} in terms of $\tau_1$ and $\eta$ therefore gives
\begin{align}
S(-(1 +12\eta + 6\eta^2)\tau_1,\eta) = S(\tau_1,\eta)\Big(1 -&\alpha \sqrt{2 + 12\eta + 6 \eta^2} + (1 + \sqrt{6}(1+\eta))S(\tau_1,\eta)\Big)\nonumber \\
\times \Big(&1+\alpha \sqrt{2 + 12\eta + 6 \eta^2}- (1 + \sqrt{6}(1+\eta))^2S(\tau_1,\eta)^2 \nonumber \\
&- S(\tau_1,\eta)(1 + \sqrt{6}(1+\eta))\Big(1 - \alpha \sqrt{2 + 12\eta + 6 \eta^2}\Big)\Big).\label{eq:4-pergov-tau-eta}
\end{align}
Analogously to the analysis of the 2-periodic case in Appendix \ref{A:logistic-transseries}, the next step is to expand this expression as a power series in $\eta$, and apply the series expression for $S(\tau_1,\eta)$ given in \eqref{eq:4per-general-ansatz}. Matching powers of $\eta^{j/2}$ for $j = 1,\ldots 4$ produces a system of four equations -- two of these equations are algebraic, and two are nonlinear ordinary differential equations in $\tau_1$. We omit the details of this step here, as it requires only algebraic manipulations, and the intermediate mathematical expressions are quite lengthy. These four equations may be simplified using the symmetry relations in \eqref{eq:4per-Theta12}--\eqref{eq:4per-Theta34}, resulting in the following system of equations
\begin{align}
\label{eq:Theta40}\Theta_{4,0}(\tau_1) &= 2 a \Theta_{1,0}(\tau_1)\Theta_{3,0}(\tau_1),\\
\label{eq:Theta20}\Theta_{2,0}(\tau_1) &= a \Theta_{1,0}(\tau_1)^2 + a\Theta_{3,0}(\tau_1)^3,\\
\label{eq:Theta10}\tau_1 \Theta_{1,0}'(\tau_1) &= \Theta_{1,0}(\tau_1) - b (\Theta_{3,0}(\tau_1)^3 + 3 \Theta_{1,0}(\tau_1)^2\Theta_{3,0}(\tau_1)),\\
\label{eq:Theta30}\tau_1 \Theta_{3,0}'(\tau_1) &= \Theta_{3,0}(\tau_1) - b (\Theta_{1,0}(\tau_1)^3 + 3 \Theta_{3,0}(\tau_1)^2\Theta_{1,0}(\tau_1)),
\end{align}
where
\begin{align}
    a  =  \frac{1}{2}\left(2+2\sqrt{6}-3 \alpha (\sqrt{2}+2\sqrt{3})\right), \qquad  b  =  \frac{5}{6}\left(14+4\sqrt{6}- \alpha (7\sqrt{2}+4\sqrt{3})\right).
\end{align}
By substituting the power series \eqref{eq:4per-general-ansatz} into the governing equation \eqref{eq:4-pergov-tau-eta}, it can be seen at leading order as $\eta \rightarrow 0$ and $\tau_1 \rightarrow 0$ that $S_{1,0} = -S_{3,0}$, providing one initial condition for the system \eqref{eq:Theta40}--\eqref{eq:Theta30}. The second initial condition may be chosen arbitrarily, as this choice may be absorbed into the expression for $\sigma_1$, in the same manner as the constant $C$ in \eqref{eq:OmegaSys}. For algebraic convenience, and without loss of generality, we select $S_{1,0} = 1$. These conditions are sufficient to uniquely solve \eqref{eq:Theta40}--\eqref{eq:Theta30}. The solution to this system is given by 
\begin{align}\label{eq:T12}
    \Theta_{1,0}(\tau_1)  = & \frac{\alpha \tau_1}{\sqrt{2-2b^2\,\tau_1^4}}\sqrt{1+\sqrt{1-b^2\tau_1^4}} ,\qquad & \Theta_{2,0}(\tau_1)  = & -\frac{a\tau_1^2}{1-b^2\tau_1^4},\\
     \Theta_{3,0}(\tau_1)  = & -\frac{\alpha b \tau_1^3}{\sqrt{2-2b^2\,\tau_1^4}}\left(\sqrt{1+\sqrt{1-b^2\tau_1^4}} \right)^{-1}, \qquad & \Theta_{4,0}(\tau_1)  = & \frac{ab \tau_1^4}{1-b^2\tau_1^4}.\label{eq:T34}
\end{align}
In principle, we can match the expansion of \eqref{eq:4-pergov-tau-eta} at higher powers of $\eta$ in order to obtain $\Theta_{j,k}$ for $j = 1,\ldots 4$ with $k > 0$. For the purposes of this example, however, the first four terms of the series will produce a useful approximation for the solution behaviour.

The final step is to determine the behaviour of the transseries parameter $\sigma_1$. This is slightly more complicated than in the 2-periodic problem, as we must incorporate the behaviour of $\hat{R}(x,\eps)$ into the calculations. We include the details of this process in Appendix \ref{A.2}, where we show that
\begin{equation}\label{per4-sigma1}
\sigma_1 = - \frac{1}{50}\left(3\sqrt{2} - 16\sqrt{3} - 7 \sqrt{6} + 12\right) + \frac{\eta}{500}\left(297\sqrt{2} - 709\sqrt{3} - 189\sqrt{2} + 399\right) + \mathcal{O}(\eta^2),
\end{equation}
We have now determined enough transseries terms to accurately approximate the solution behaviour in the 4-periodic regime.

\subsection{Error comparison}

\begin{figure}[tb]
\centering
\subfloat[Exact solution (red circles) and 2-periodic approximation (blue dots) for $\eps = 0.05$.]{
\includegraphics[width=0.425\textwidth]{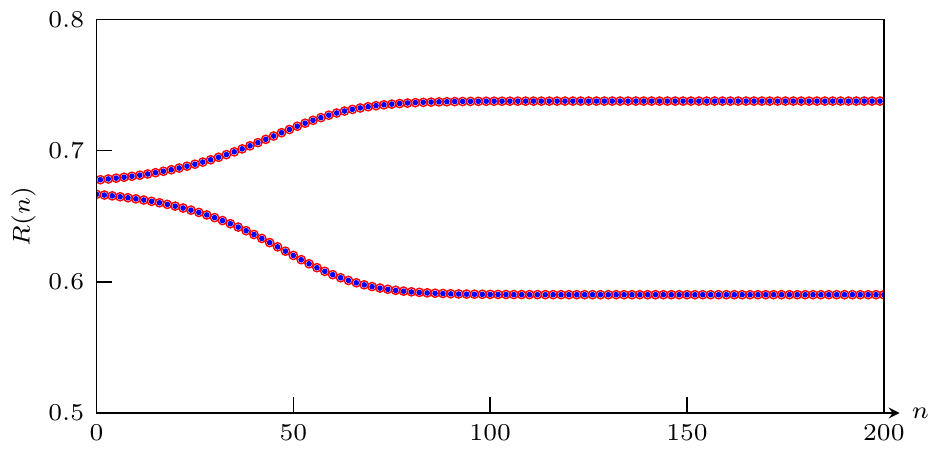}
}
\subfloat[Exact solution (red circles) and 4-periodic approximation (blue dots) for $\eps = 0.51$.]{

\includegraphics[width=0.425\textwidth]{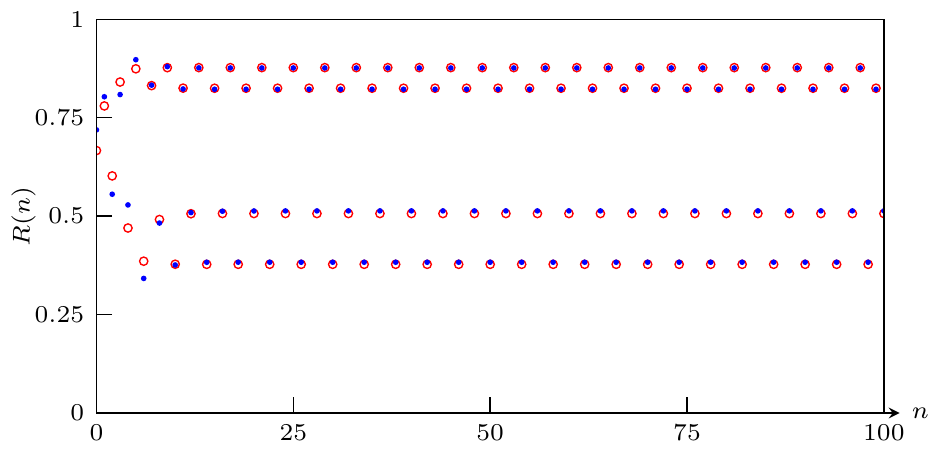}
}

\subfloat[2-periodic approximation error for $\eps = 0.05$.]{

\includegraphics[width=0.425\textwidth]{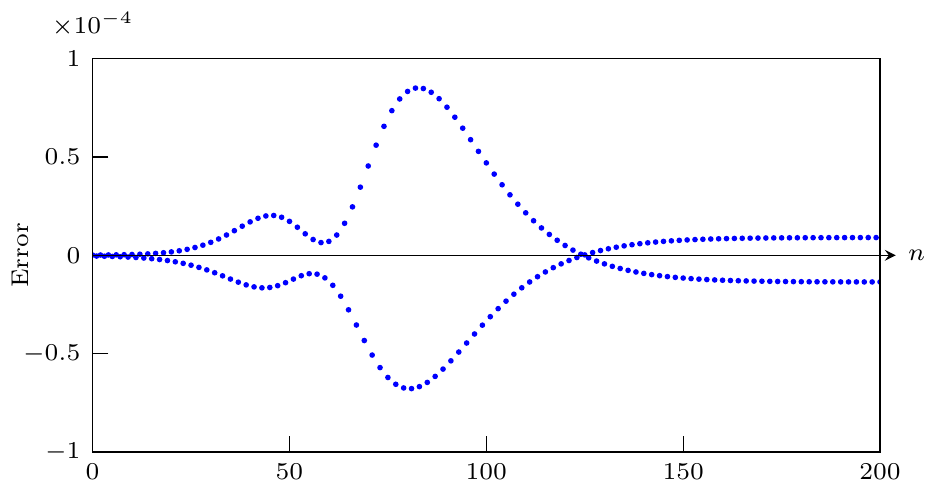}
}
\subfloat[4-periodic approximation error for $\eps = 0.51$]{

\includegraphics[width=0.425\textwidth]{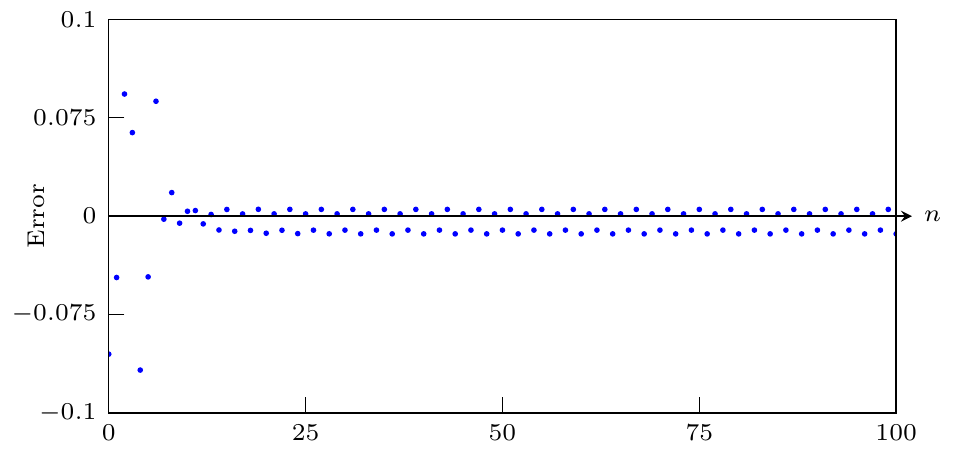}
}

\caption{The plot in (a) compares the 2-periodic approximation $R_{2,\mathrm{app}}$, from \eqref{eq:rapp2}, against the exact solution for $\eps = 0.05$. The plot in (b) compares the 4-periodic approximation $R_{4,\mathrm{app}}$, from \eqref{eq:rapp4}, against the exact solution for $\eps = 0.51$, or $\sqrt{6}\eta \approx 0.0605$. The approximation errors, given by the difference between the exact solution $R(x)$ and the approximations are shown in (c) and (d). The 2-periodic approximation has maximum error in the region just before reaching the 2-periodic steady solution. The 4-periodic approximation has maximum error in the initial region; this is to be expected, as the initial condition was obtained directly from the 2-periodic solution, and is not expected to be highly accurate in the 4-periodic regime.}\label{fig:static2-example}
\end{figure}

As a consequence of the preceding analysis, we are able to derive an approximation for the solution to the logistic equation in the 2-periodic and 4-periodic parameter regimes, which we denote as $R_{2,\mathrm{app}}(x)$ and $R_{4,\mathrm{app}}(x)$ respectively. Combining \eqref{eq:logistic-tau-constant}, \eqref{eq:logistic-R-newdef}, \eqref{eq:static-omega0}, \eqref{eq:static-omega1}, and \eqref{eq:per2-sigma1}, we find that in the 2-periodic parameter regime
\begin{equation}\label{eq:rapp2}
R(x) \approx R_{2,\mathrm{app}}(x) = \frac{2+\eps}{3+\eps} + \eps^{1/2}\Omega_{o,0}({\tau_0}) + \eps \Omega_{e,0}({\tau_0})+\eps^{3/2}\Omega_{o,1}({\tau_0}) + \eps^2 \Omega_{e,1}({\tau_0}),
\end{equation}
where $\tau_0$ and $\sigma_0$ are approximated as
\begin{equation}
{\tau_0} = \sigma_0 \eps^{1/2}\e^{-x(\pi \ii   +  \log(1 + \eps))/\eps},\qquad \sigma_0 \approx -\frac{1}{9} + \frac{4\eps}{81} - \frac{19\eps^2}{648}.
\end{equation}
A comparison of the exact solution against the approximation is shown for $\eps = 0.05$ in Figure \ref{fig:static2-example}(a). The exact solution is shown as red circles, while the approximation is shown as blue dots. The two curves are visually indistinguishable. The approximation error is shown in Figure \ref{fig:static2-example}(c). It is clear from this figure that the error has a peak at the end of the transition region, just before the solution settles into the stable 2-periodic behaviour. 

In the 4-periodic parameter regime $\eps>-2+\sqrt{6}$, the approximated transseries is given combining the expressions in \eqref{eq:4per-tau1}, \eqref{eq:T12}--\eqref{per4-sigma1}, and the previous approximation \eqref{eq:rapp2}, to give
\begin{align}\label{eq:rapp4}
R(x) \approx R_{4,\mathrm{app}}(x) = R_{2,\mathrm{app}}(x) + \sqrt{\eta}(\Theta_{1,0}(\tau_1) +& \Theta_{3,0}(\tau_1)) + \eta (\Theta_{2,0}(\tau_1) + \Theta_{4,0}(\tau_1)) ,
\end{align}
where $\tau_1$ and $\sigma_1$ are approximated as
\begin{align}
{\tau_1} &= \sigma_1 \eps^{1/2}\e^{x(\log(1-\eps(4+\eps))/2 + \pi \ii)/\eps},\\
 \sigma_1 &\approx - \frac{1}{50}\left(3\sqrt{2} - 16\sqrt{3} - 7 \sqrt{6} + 12\right)+ \frac{\eta}{500}\left(297\sqrt{2} - 709\sqrt{3} - 189\sqrt{2} + 399\right)  .
\end{align}
Note that we do not include the term containing $g(x,\eps)$ in $B(x,\eps)$ from \eqref{eq:4per-A1}. This term disappears for integer values of $n$, and therefore can be omitted at this stage without altering the approximation.

\begin{figure}[tb]
\centering
\subfloat[Error for 2-Periodic Approximation]{

\includegraphics[height=0.4\textwidth]{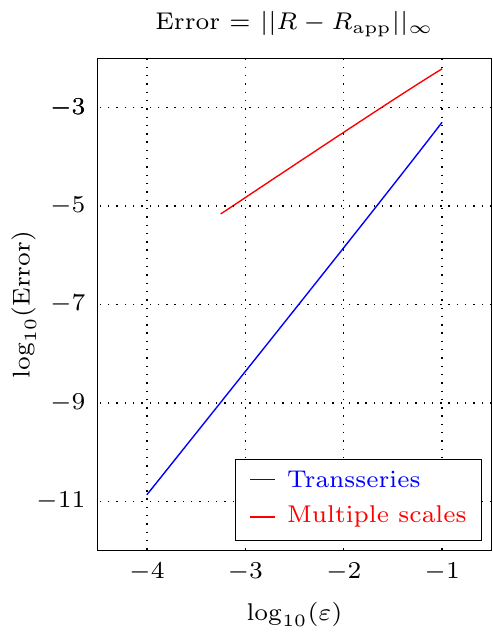}
}
\subfloat[Error for 4-Periodic Approximation]{

\includegraphics[height=0.4\textwidth]{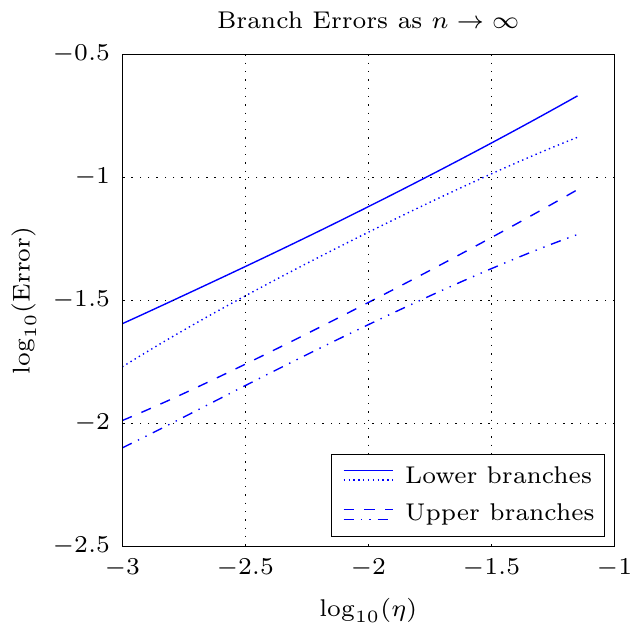}
}
\caption{The plot in (a) shows the resummed transseries approximation error in blue, corresponding to the maximum difference between the approximated and exact value. This measure of the error was chosen to be consistent with the error measure provided in \cite{Hall2016}; this error is shown as a red curve. Due to the ease with which the transseries method captures higher-order exponential behaviour, which plays an important role in the transition region between non-periodic and 2-periodic behaviour, it outperforms the multiple scales approximation. In (b), we show the error of the four solution branches as $n \rightarrow \infty$, approximated by taking $|R - R_{4,\mathrm{app}}|$ on each of the four branches for a value of $n$ sufficiently large that the error is not visibly changing. For each of the four branches, the error decreases as $\eta \rightarrow 0$, as would be expected.}\label{fig:static-approx}
\end{figure}

In Figure \ref{fig:static-approx}(a), we show the approximation error for a range of values of $\eps$, where the error is measured as the maximum difference between the exact solution and the transseries approximation, shown as a blue curve. This error measure was chosen to allow for direct comparison with equivalent results from \cite{Hall2016}, which are shown as a red curve. The transseries approximation is more accurate than the multiple scales approximation in this parameter regime, and the error decays faster in the limit that $\eps \rightarrow 0$. The reason for this behaviour is that the transseries approach allowed for higher-order exponential corrections to be easily computed and retained. The maximum approximation error occurs at the end of the transition region between non-periodic and 2-periodic behaviour, where the exponential contributions contribute significantly to the solution behaviour. Computing these exponential corrections using multiple scales methods would be an algebraically significantly more demanding task, requiring matched asymptotic expansions to be applied at higher orders of the expansion.

\subsection{8-periodic solution}\label{S:8-per}

We may continue this process to understand the emergence of the next period doubling bifurcation. While we will not include a full explicit, algebraic analysis here, we will show that the exponential factor can be used to identify the appearance of 8-periodic stable solutions as $\eps$ is increased further.

The method from  Section \ref{s:static-4per} can be applied again in order to obtain approximations for solutions with even higher periodicity. We can now write the next term in the transseries such that $R(x,\eps) = \hat{R}(x,\eps) + S(x,\eps) + T(x,\eps)$. The quantity $T(x,\eps)$ is defined in terms of a new transseries parameter $\sigma_2$ to be 
\begin{equation}
T(x,\eps; \sigma_2) = \sum_{m=1}^{\infty} \sigma_2^m \e^{-mF(x,\eps)/\eps} T_m(\eps).\label{eq:logistic-transseries-constant8}
\end{equation}
The transseries terms $\hat{R} + S$ capture the 4-periodic solution behaviour, and therefore must tend to the 4-periodic solution in the limit that $\tau_0$ and $\tau_1$ become large. We denote this solution as $R_4(\eps)$. Hence, we apply the expression $R(x,\eps) = R_4(\eps) + T(x,\eps)$ to the governing equation \eqref{eq:logistic-finite-diff-constant} and find an expression for the exponential weights, in similar fashion to the process for obtaining \eqref{eq:logistic-A} or \eqref{eq:static-bd}. 

\begin{figure}
\centering

\includegraphics{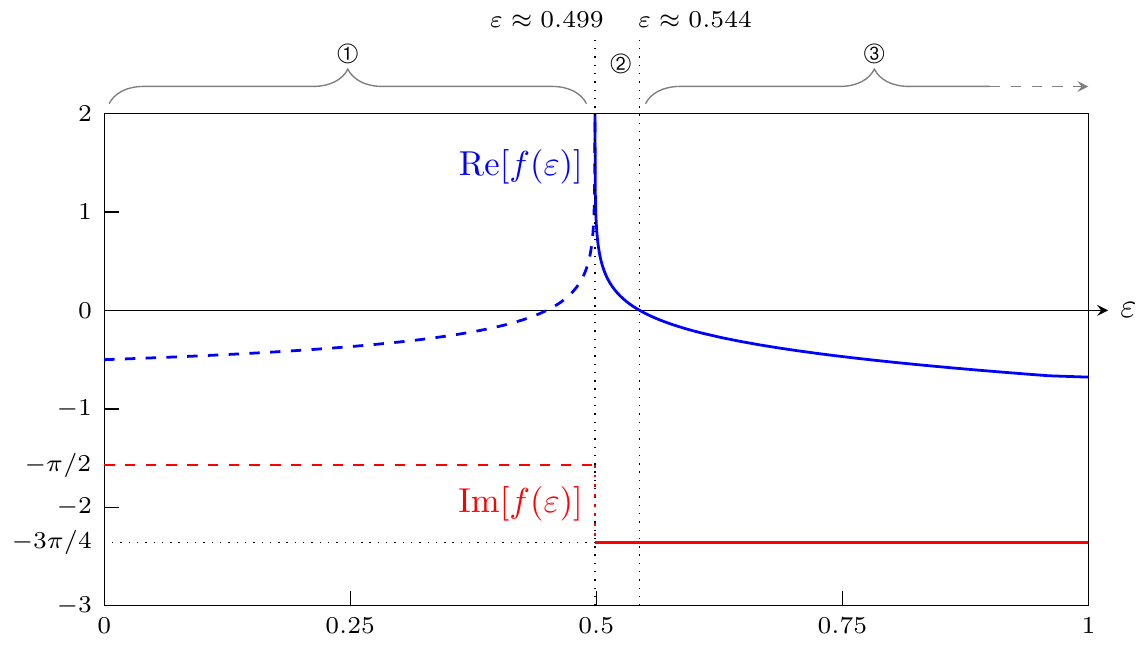}

\caption{This figure shows the real and imaginary parts of $f(\eps)$, where $F(x,\eps) = f(\eps) x + \eps g(x,\eps)$. The behaviour of the transseries depends on both the real and imaginary part of this quantity, in the same fashion as Figure \ref{fig:A1-log}. The exponential contribution is not present in the transseries in parameter regime \ding{192}. In regimes \ding{193} and \ding{194} the contribution is present, and must be 8-periodic, due to the value of $\mathrm{Im}[f(\eps)]$ in these regimes. In regime \ding{193}, the 8-periodic contribution is small, due to the positive sign of $\mathrm{Re}[f(\eps)]$, and this contribution becomes significant in regime \ding{194}, as the sign of $\mathrm{Re}[f(\eps)]$ becomes negative.}\label{fig:A2-log}
\end{figure}

The exponential weights may again be written in the form $F(x,\eps) = f(\eps) x + \eps g(x,\eps)$, where $g$ disappears on $n \in \mathbb{Z}$. The behaviour of $f(\eps)$ is illustrated in Figure \ref{fig:A2-log}. A very similar set of inferences may be drawn from this image as for Figure \ref{fig:A1-log}. In region \ding{192}, the 8-periodic behaviour does not contribute to the solution, as discussed for the 4-periodic case in Section \ref{S:where-4-periodic-behaviour-matters}. This 8-periodic contribution appears in the transseries as $\eps$ moves into region \ding{193}. 
In this range of $\eps$, there are 8-periodic contributions to the solution, but they are smaller than the 4-periodic solution contribution, as the exponential term is relatively small compared to those in $S(x,\eps)$, decaying exponentially as $x\to\infty$. Finally, in region \ding{194}, the 8-periodic solution grows exponentially, and the behaviour of $T(x,\eps)$ dominates the solution behaviour. 

It is therefore clear that we can explain the onset of these higher periodicity solutions by explicitly studying the exponential weights of the transseries solution; while the algebraic complexity of the process increases after each doubling, the steps for identifying this behaviour remain essentially the same. The resummed transseries therefore provides a systematic approach to studying bifurcations even for larger values of the bifurcation parameter, where classical asymptotic methods typically fail.

\section{Dynamical Logistic Equation}\label{s:dynamic}

In the previous section, we studied the classical logistic equation, and showed that the higher periodicity solutions may be obtained directly using a transseries approach. In this section, we consider a more complicated variant of this problem, known as the slowly-varying logistic equation. 
\begin{equation}\label{eq:dynamic-logistic-original}
y(n+1) = (\lambda_0 + \eps n) y(n) [1-y(n)],\qquad 0 < y(0) < 1,
\end{equation}
with $\eps > 0$. The bifurcation parameter is given by $\lambda = \lambda_0 + \eps n$, and in this case, it changes slowly over time. In previous studies \cite{Fruchard1996, Fruchard2009, ElRabih2003}, this has been shown as an example of a ``canard'' solution, in which the behaviour appears to remain near the unstable solution for an extended period of time, before rapidly jumping to approach the stable solution with higher periodicity. As $n$ increases, this parameter will pass through values across which the solution stability is known to change. When $\lambda_0 + \eps n = 3$, the 1-periodic equilibrium becomes unstable, and the 2-periodic equilibrium becomes stable. As $n$ increases further, eventually $\lambda$ exceeds $1 + \sqrt{6}$, and the 2-periodic equilibrium becomes unstable, with the 4-periodic equilibrium becoming stable. This process continues until the bifurcation parameter becomes sufficiently large that the solution becomes chaotic. For the problem studied here, we will set $\lambda_0 = 3$ and $y(0) = 2/3$.

In \cite{Hall2016}, it was shown that a discrete multiple scales approach can be used to describe this behaviour asymptotically. This approach required balancing several different timescales, and using asymptotic matching to connect the solutions in each different asymptotic region.

In this section, we will show that this process can be described using a transseries approach, with the resulting expansion to be valid even as the solution behaviour changes dramatically, and increases in periodicity.  We will now show that transseries  provide a systematic and generally more accurate approach than the multiple scales procedure of \cite{Hall2016} in describing the solution behaviour as it transitions from an unstable to a stable manifold; this will demonstrate that transseries expansions can be used to effectively capture canard behaviour in discrete systems.

We will show the first stability transition in detail. We will subsequently provide an outline of how this method can be extended to describe the second transition, together with some results; however, the algebraic manipulations for this process are quite involved, and the precise details will be omitted.

\subsection{2-periodic solution}

\subsubsection{Transseries ansatz}\label{s:dyn-2per-ansatz}

The difference from \eqref{eq:logistic-finite-diff-constant} above is that in the prefactor of the r.h.s. the perturbative parameter $\eps$ is replaced now by $x$. We again begin by applying a multiple scales ansatz, and expanding as a transseries in a continuous variable $x$. Setting $x = \eps n$ and $R(x) = y(n)$ gives
\begin{equation}
R(x+\eps) = (3+x)R(x)[1-R(x)],\qquad R(0) = \tfrac{2}{3}.\label{eq:dynamic}
\end{equation}
We again formulate an ansatz for the solution in terms of $\eps$ and a transseries parameter $\sigma_0$. The ansatz is identical to that given in \eqref{eq:logistic-transseries-constant}, but has been included below for convenience:
\begin{equation}
R(x,\eps; \sigma_0) = \sum_{m=0}^{\infty} \sigma_0^m \e^{-mA(x)/\eps} R_m(x,\eps),\qquad R_m(x,\varepsilon) \simeq \varepsilon^{\beta_m}\sum_{k=0}^{\infty}\varepsilon^k R_{m,k}(x),\label{eq:logistic-transseries-constant2}
\end{equation}
where $\beta_m$ will again be chosen such that $R_{m,0}$ takes nonzero value.
It is straightforward to compute the first few terms of the algebraic portion of the series expression, corresponding to $m=0$ in \eqref{eq:logistic-transseries-constant2}, which gives a power series expression for the non-periodic manifold. The recursion relation is given obtained by expanding $R(x+\eps)$ using a power series in $\eps$, and matching powers of $\eps$ in the resultant expression. This process gives
\begin{align}
R_{0,0}(x) & =\frac{2+x}{3+x},\label{eq:logistic-pert-coeffs}\\
R_{0,k}(x) & =-\frac{1}{(2+x)}\left[\sum_{n=1}^{k}\frac{1}{n!}\,R_{0,k-n}^{(n)}(x)+(3+x)\sum_{n=1}^{k-1}R_{0,n}(x)R_{0,k-n}(x)\right],\qquad k\ge1.
\end{align}
The first few iterations of this recurrence relation give
\begin{equation}
R_{0,1}(x) = -\frac{1}{(x+2)(x+3)^{2}},\qquad R_{0,2}(x) = \frac{x^{2}+x-4}{{(x+2)^{3}(x+3)^{3}}},\qquad R_{0,3}(x) =  -\frac{x^{4}-2x^{3}-28x^{2}-33x+24}{(x+2)^{5}(x+3)^{4}}.
\end{equation}
This process may be continued indefinitely in order to continue calculating terms in the power series for the non-periodic manifold. This process will not, however, capture the transition to the 2-periodic manifold. In order to obtain an approximation for this behaviour, we are required to consider terms in the ansatz \eqref{eq:logistic-transseries-constant2} with $m \neq 0$. Continuing to the next order in $\sigma_0$, we find that
\begin{equation}
\e^{-[A(x+\eps) - A(x)]/\eps} R_1(x + \eps,\eps) = (3+x) R_1(x,\eps)[1-2R_0(x,\eps)].
\end{equation}
As before, the argument of the exponential may be determined by expanding $R_1$ as a power series in $\eps$, as well as expanding $A(x + \varepsilon) = A(x) + \varepsilon A'(x) + \cdots$. At leading order in $\eps$, this gives the differential equation
\begin{equation}
\mathrm{e}^{-A'(x)}=-\left(x+1\right)=\mathrm{e}^{-\left(2p+1\right)\pi\mathrm{i}+\log\left(x+1\right)},\qquad p\in\mathbb{Z}.
\end{equation}
Hence, we obtain
\begin{equation}
A(x)=(2p+1)\pi\mathrm{i}\,x+x-(x+1)\log(x+1),\label{eq:dyn-logistic-action-general}
\end{equation}
where we follow the same reasoning as the analysis used to determine \eqref{eq:logistic-A}, and absorb the constant into the series parameter. We may again set $p=0$; this choice will have no effect on the behaviour of the solution for integer values of $n$.

Once $A(x)$ has been determined, it is possible to obtain a recurrence relation for $R_m(x)$ by applying the first ansatz expression in \eqref{eq:logistic-transseries-constant2} to the governing equation \eqref{eq:dynamic}, and matching powers of the transseries parameter $\sigma_0$. This gives
\begin{align}
\nonumber(-1)^{m}(1+x+\varepsilon)^{m}\mathrm{e}^{m((1+x)\log(1+{\varepsilon}/{(1+x)})/{\varepsilon}-1)}&R_{m}(x+\eps,\eps)\\
=(3+x)R_{m}(x,\eps)&[1-2R_{0}(x,\eps)]-(3+x)\sum_{n=1}^{m-1}R_{n}(x,\eps)R_{m-n}(x,\eps).\label{eq:logistic-higher-sectors-rec}
\end{align}
It is now possible to apply the second part of the ansatz in \eqref{eq:logistic-transseries-constant2} and to match powers of $\eps$ in this expression. By direct substitution, we find that $\beta_m = m$ gives the result that $R_{m,0}$ is nonzero. By subsequently matching terms which are $\mathcal{O}(\eps)$ in the small $\eps$ limit, it is possible to generate an equation for $R_{1,0}$ and a recurrence relation for $R_{m,0}$ for $m\geq 2$. We find that
\begin{align}
(x+1)R'_{1,0}(x) = -\left(\frac{2}{x+2} - \frac{2}{x+3} + \frac{1}{2}\right)R_{1,0}, \label{eq:dyn_R10}
\end{align}
The initial condition in \eqref{eq:dyn_R10} may be chosen arbitrarily, as this choice may be absorbed into the transseries parameter. Choosing $R_{1,0}(0) = 1$ gives
\begin{equation}
R_{1,0}(x) = \frac{3(x+2)^2}{4(x+1)^{3/2}(x+3)}.\label{eq:dyn-R10}
\end{equation}
The recurrence relation for subsequent terms is given by
\begin{equation}
\left[(-1)^m(1+x)^m + (1+x)\right]R_{m,0}(x) = -(3+x)\sum_{n=1}^{m-1}R_{n,0}(x)R_{m-n,0}(x),\qquad m \geq 2.
\end{equation}
Continuing to match higher powers of $\eps$ allows for the direct computation of terms further terms such as $R_{m,k}$, obtained by matching terms which are $\mathcal{O}(\eps^k)$ in the small $\eps$ limit. The direct computation of further terms is not required for the present analysis.

\subsubsection{Computing terms in the resummed transseries}

Motivated by the analysis of the static system, and in particular the form of \eqref{eq:2per-general-ansatz}, we switch the order of summation in the transseries \eqref{eq:logistic-transseries-constant2}, writing it as
\begin{equation}
R(x,\eps;\sigma_0) \simeq \sum_{k=0}^{\infty}\eps^k\sum_{m=0}^{\infty}\left(\sigma_0 \eps \e^{-A(x)/\eps}\right)^m R_{m,k}.\label{eq:dyn-summation-switched}
\end{equation}
A main difference from the static one is that the expansion in powers of $eps$ is asymptotic, while the sum over the exponentials ($m\ge0$) is convergent. Thus \eqref{eq:dyn-summation-switched} is a formal expansion. 

As for the static system, we define a new series parameter $\tau_0$, and new quantities $\Omega_k(\tau_0)$ such that
\begin{equation}
\tau_0 = \sigma_0 \eps \e^{-A(x)/\eps},\qquad \Omega_k(\tau_0) = \sum_{m=0}^{\infty}\tau_0^m R_{m,k}.\label{eq:static-tau}
\end{equation}
It will be helpful later to note that
\begin{equation}
\tau_0(x+\eps) = \e^{-(A(x+\eps) - A(x))/\eps}\tau_0(x) = \tau_0(x)\left[\e^{-A'(x)} + \mathcal{O}(\eps)\right] \quad \mathrm{as} \quad \eps \rightarrow 0.\label{eq:dyn-tau-plus-eps}
\end{equation}
The transseries expression in \eqref{eq:dyn-summation-switched} is now given by
\begin{equation}
R(\tau_0,\eps) \simeq \sum_{k=0}^{\infty}\eps^k \Omega_k(\tau_0).
\end{equation}
We can now apply this expression to \eqref{eq:dynamic} and match orders of $\eps$. At leading order, we find that
\begin{equation}
\Omega_0\left(\e^{-A'(x)}\tau_0\right) = (3+x)\Omega_0(\tau_0)(1-\Omega_0(\tau_0)),\label{eq:dyn-static-0}
\end{equation}
where \eqref{eq:dyn-tau-plus-eps} was used to obtain the leading-order on the left-hand side. At this stage, we could mechanically obtain the function $\Omega_0$ as a Taylor series in $\tau_0$, which is convergent, with some finite radius of convergence. It happens, however, that there exists a particularly convenient variable transformation that converts the right-hand side from a dilation to a translation. If we define a new variable $y$ such that $y = -x\log(\tau_0)/A'(x)$, the expression in \eqref{eq:dyn-static-0} becomes
\begin{equation}
\Omega_0(y + x) = (3+x)\Omega_0(y)(1-\Omega_0(y)).\label{eq:MAGIC}
\end{equation}
This expression has the same form as the static logistic map equation, given in \eqref{eq:logistic-finite-diff-constant}, with $x$ in place of $\eps$. Furthermore, since $x = \eps n$, it is valid to apply the asymptotic solution derived for this expression in Section \ref{S:static-2-per}. As we are interested in capturing the first transition, across which the solution switches from having no periodic component to having a 2-periodic component, we can directly apply the transseries expression for the 2-periodic solution given in \eqref{eq:logistic-R-newdef}. 

In order to take into account the form of \eqref{eq:MAGIC}, we must replace $\eps$ and $x$ with $x$ and $y$ respectively in \eqref{eq:logistic-R-newdef}. We must also replace the $\tau_0$ in this expression with a new transseries parameter $\overline{\tau}_0$, in which $\eps$ and $x$ are again replaced with $x$ and $y$ respectively. This gives 
\begin{equation}\label{eq:overline-tau0}
    \overline{\tau}_0(y,x) = \overline{\sigma}_0\sqrt{x}\e^{-\ii\pi y/x + y \log(1+x)/x} =  \overline{\sigma}_0\sqrt{x}\e^{(\ii\pi  + \log(1+x))\log(\tau_0)/A'(x)} =  \overline{\sigma}_0\sqrt{x} \tau_0,
\end{equation}
where $\overline{\sigma}_0$ is a new transseries parameter that remains to be determined. Making the appropriate substitutions in \eqref{eq:logistic-R-newdef} now gives the form of $\Omega_0(y)$ as
\begin{equation}
\Omega_0(y) = \frac{2+x}{3+x} + \sqrt{x} \sum_{k=0}^{\infty}x^{k} \Omega_{o,k}(\overline{\tau}_0) + x \sum_{k=0}^{\infty} x^k \Omega_{e,k}(\overline{\tau}_0),\label{eq:dyn-new-omega-0}
\end{equation}
where $\Omega_{o,k}$ and $\Omega_{e,k}$ are defined in \eqref{eq:logistic-tau-constant}, and $\Omega_{o,k}$ and $\Omega_{e,k}$ for $k=0$ and $k=1$ are given explicitly in \eqref{eq:static-omega0} and \eqref{eq:static-omega1} respectively.

In a typical problem of this form, $\overline{\sigma}_0$ would be determined using the fact that $\Omega_0 = 2/3$ at $y = 0$; however, this is enforced by the transformation $y = -x \log(\tau_0)/A'(x)$, which forces $x$ to be zero if $y=0$. Consequently, the initial condition cannot be used to determine $\overline{\sigma}_0$. This is to be expected, as the the initial condition will instead be used to determine the original transseries parameter $\sigma_0$. 

Instead, we expand \eqref{eq:dyn-new-omega-0} as a Taylor series about $x = 0$ using the form of $\Omega_{o,0}$ and $\Omega_{e,0}$ given in \eqref{eq:static-omega0}. This gives
\begin{equation}
\Omega_0(y) = \frac{2+x}{3+x} + \overline{\sigma}_0 x \tau_0 - \frac{3}{2}(\overline{\sigma}_0 x \tau_0)^2 + \ldots,
\end{equation}
where the omitted terms are $\mathcal{O}(x^3 \tau_0^2)$. We may now match powers of $\tau_0$ with \eqref{eq:dyn-summation-switched} to determine that $x \overline{\sigma}_0 = R_{1,0}$, which was explicitly calculated in \eqref{eq:dyn-R10}. We therefore find that 
\begin{equation}\label{eq:overline-sigma}
\overline{\sigma}_0 = \frac{3(x+2)^2}{4x(x+1)^{3/2}(x+3)}.
\end{equation}

\subsubsection{Initial Condition}

We have now explicitly calculated all of the required quantities for the transseries approximation except for $\sigma_0$, which must be determined from the initial condition at $x=0$. At $x = 0$, it follows that $\tau_0 = \sigma_0 \eps$. Consequently, the initial condition is given by $R(\tau_0=\sigma_0\eps,x=0) = 2/3$, which we apply to the first expression in \eqref{eq:logistic-transseries-constant2}. We then express $\sigma_0$ as a power series in $\eps$, where the series terms are functions of $R_{m,k}$ for various values of $m$ and $k$, giving
\begin{equation}
\sum_{m=0}^{\infty}\sigma_0^m R_m(0,\eps),\quad \mathrm{where} \quad \sigma_0 = \sum_{j = 0}^{\infty}\eps^{j} \sigma_{0,j}.
\end{equation}
Using the second expression from \eqref{eq:logistic-transseries-constant2} and matching powers of $\eps$ allows us to compute $\sigma_{0,j}$. We have obtained enough $R_{m,k}$ terms to solve for $\sigma_{0,0}$, giving
\begin{equation}\label{eq:dyn-initial-cond}
\sigma_0 = -R_{0,1} + \mathcal{O}(\eps) = \frac{1}{18} + \mathcal{O}(\eps).
\end{equation}
Computing subsequent series terms for $\sigma_0$ requires values of $R_{m,k}$ that are not presented in this study, as even this first order approximation is sufficiently accurate as we now show.

\subsection{Error comparison}

\begin{figure}[tb]
\centering
\subfloat[Exact solution (red circles) and transseries approximation (blue dots) for $\eps = 0.001$.]{

\includegraphics[width=0.8\textwidth]{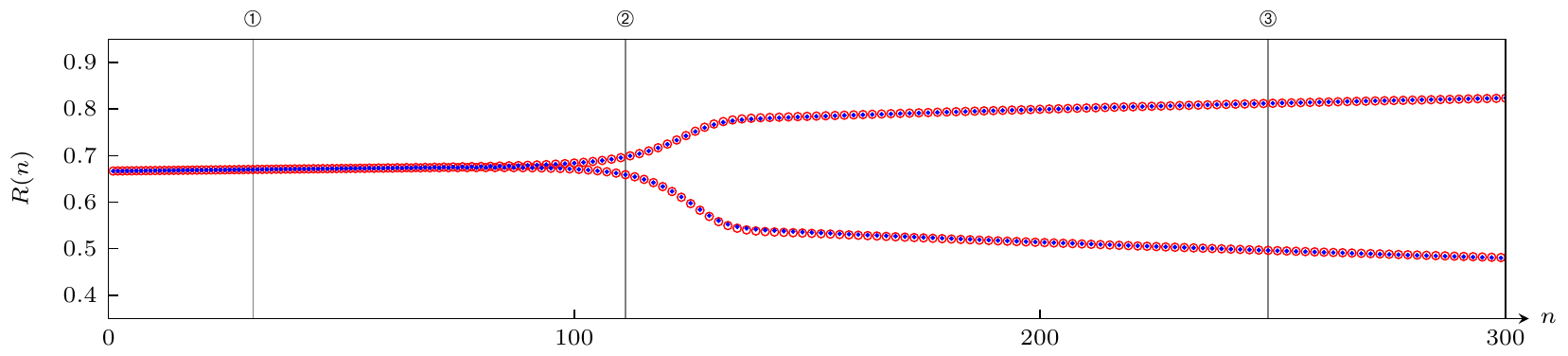}
}

\subfloat[Approximation error for $\eps = 0.001$.]{

\includegraphics[width=0.8\textwidth]{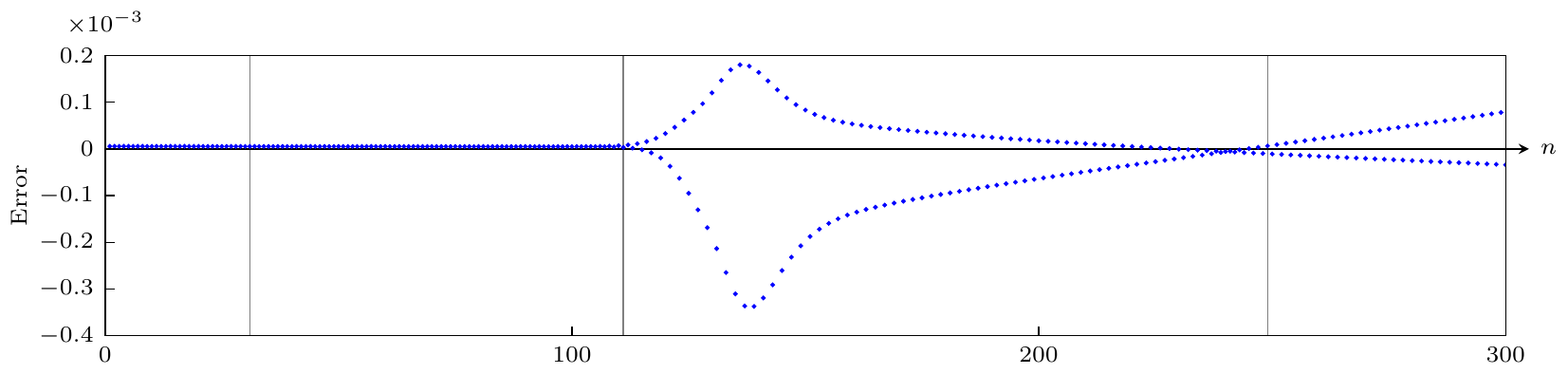}
}

\caption{The plot in (a) shows the approximation \eqref{eq:dynamic-2per-approx} and exact solution of \eqref{eq:dynamic-logistic-original} for $\eps = 0.001$. The difference between these is shown in (b). The points labeled \ding{192}, \ding{193} and \ding{194} will be referenced below in Figure \ref{fig:dyn-error}. We see that the error reaches a maximum at the start of the 2-periodic regime. It then decreases, although will eventually increase as $n$ grows, due to the increasing influence of the 4-periodic region which was not computed. Note that the small error at the point labelled \ding{194} corresponds to the value where the approximation crosses the exact solution. This occurs at some point in the 2-periodic region for any choice of $\eps$, and therefore does not signify a special parameter choice. It is an artifact of the error calculation.}\label{fig:dyn-approx}
\end{figure}

As a consequence of the preceding analysis, we are able to derive an approximation for the solution to the slowly-varying logistic equation \eqref{eq:dynamic}, which we denote $R_{\mathrm{app}}$. Combining \eqref{eq:logistic-tau-constant}, \eqref{eq:static-tau}, \eqref{eq:overline-tau0}, \eqref{eq:dyn-new-omega-0}, \eqref{eq:overline-sigma} and \eqref{eq:dyn-initial-cond}, we find an approximation for the transseries solution
\begin{equation}\label{eq:dynamic-2per-approx}
R(x) \approx R_{\mathrm{app}}(x) = \frac{2+x}{3+x} + x^{1/2}\Omega_{o,0}(\overline{\tau}_0) + x \Omega_{e,0}(\overline{\tau})+ x^{3/2}\Omega_{o,1}(\overline{\tau}_0) + x^2 \Omega_{e,1}(\overline{\tau}_0),
\end{equation}
where
\begin{equation}
\overline{\tau}_0 \approx \frac{\epsilon (x+2)^2 \e^{-(\pi \ii  x + x - (x+1)\log(x+1))/\eps}}{24x^{1/2} (x+1)^{3/2}(x+3)}.
\end{equation}
The most useful feature of this approximation is that it is valid before, during, and after the transition region in the slowly-varying logistic equation. We illustrate an example comparison in Figure \ref{fig:dyn-approx}(a), corresponding to $\eps = 0.001$. The approximation is shown as blue dots, and overlaid on top of the exact solution, shown as red circles. The two solutions are visually almost indistinguishable. 

The approximation error for this example is shown in Figure \ref{fig:dyn-approx}(b), calculated by $y(n) - R_{\mathrm{app}}(\eps n)$. The error reaches a peak following the transition region, at the beginning of the stable 2-periodic behaviour. The error does grow in this region as $n$ becomes large, and continues to do so until the transition to 4-periodic behaviour occurs. This behaviour is not depicted in Figure \ref{fig:dyn-approx}(b).

\begin{figure}[tb]
\centering
\subfloat[Non-periodic region at \ding{192}]{

\includegraphics[height=0.35\textwidth]{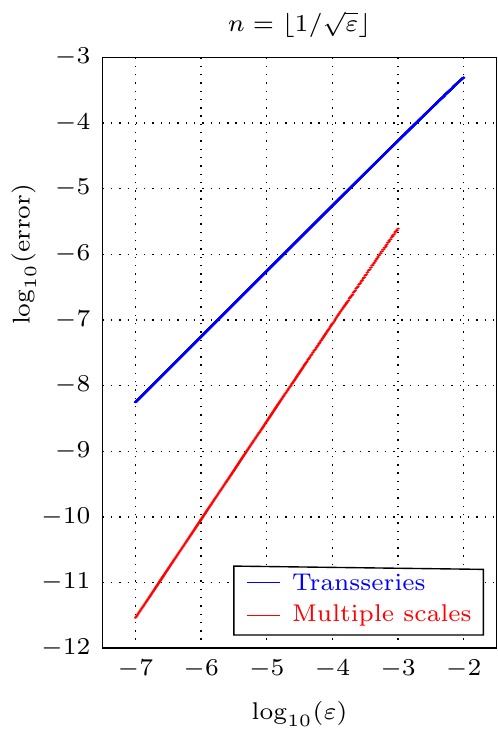}
}
\subfloat[Transitional region at \ding{193}]{

\includegraphics[height=0.35\textwidth]{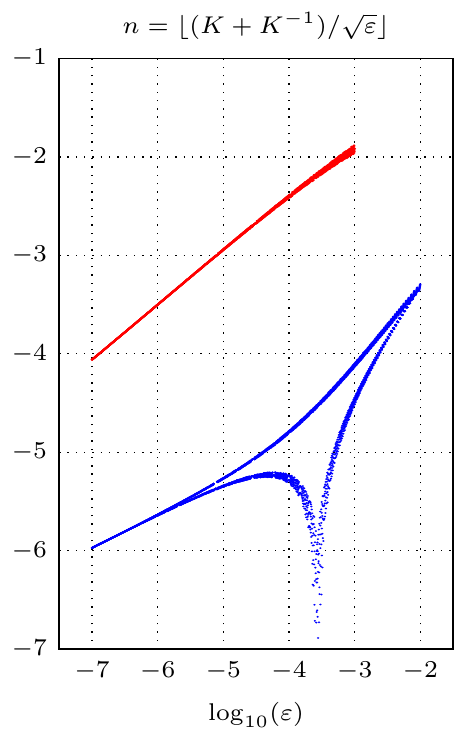}
}
\subfloat[2-periodic region at \ding{194}]{

\includegraphics[height=0.35\textwidth]{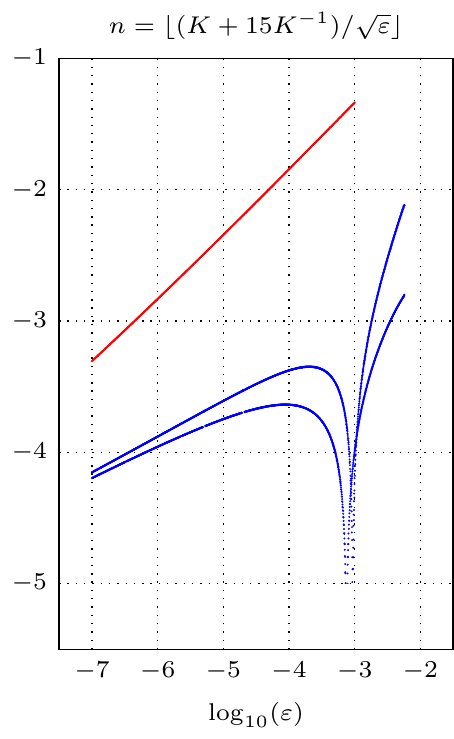}
}
\caption{This figure shows the error in the dynamic system at the three points identified in \cite{Hall2016} as belonging to the inner region, transition region, and outer region, shown as points \ding{192}, \ding{193} and \ding{194} in Figure \ref{fig:dyn-approx}. In each case, the error is shown as a blue curve. This curve becomes smaller as $\varepsilon$ decreases. The point at which the error dips is an artifact of the observation that the approximation crosses the exact solution at the identified point for this choice of $\eps$, and does not represent any significant phenomenon within the transseries approximation. The cause of this behaviour is explained in more detail within the description of Figure \ref{fig:dyn-approx}. We have chosen a similar range of small parameter to the analysis in \cite{Hall2016}, shown in red. The transseries outperforms the multiple scales method in both the transition region and the outer region, in which the exponential terms play an important role in describing the solution behaviour. These terms are more naturally captured using transseries methods, leading to an improved approximation.}\label{fig:dyn-error}
\end{figure}

In order to obtain a more complete picture of the accuracy of the transseries approximation, we determined the approximation error at three selected values of $n$. These values were tested in \cite{Hall2016} relative to other methods, to obtain representative computations of the approximation error in important parts of the solution domain. The first point is $n = \lfloor 1/\sqrt{\eps} \rfloor$. This point is found in the early non-periodic region before the transition from non-periodic behaviour to 2-periodic behaviour occurs. It is labelled \ding{192} in the example solution from Figure \ref{fig:dyn-approx}.

For comparison, we need to identify the remaining representative points used in \cite{Hall2016}, which required the computation of an intermediate quantity $K$, satisfying
\begin{equation}
K = \sqrt{\log K - \tfrac{3}{2}\log(\eps)}.
\end{equation}
This quantity was derived in \cite{Hall2016} although it has been adjusted to take into account the slightly different form for the slowly-varying logistic equation considered here. The second point falls within the transition region between the non-periodic unstable manifold and the 2-periodic stable manifold, and is given by $n = \lfloor K + K^{-1}/\sqrt{\eps} \rfloor$. This point is labelled \ding{193} in the example solution from Figure \ref{fig:dyn-approx}. Finally, we also determine the error at a point in the region where the solution has completed its transition to 2-periodic stable behaviour. This point is given by $n = \lfloor K + 15 K^{-1}/\sqrt{\eps} \rfloor$, and is labelled \ding{194} in the example solution from Figure \ref{fig:dyn-approx}. The error for each of these three points was studied in \cite{Hall2016} allowing for direct comparison between the transseries approximation and the multiple scales approximation errors.

The error for each of the three representative points over a range of $\eps$ values may be seen in Figure \ref{fig:dyn-error}, shown in blue. The error for the approximation from the multiple scales approximation in \cite{Hall2016} is shown in red for each point. In each region, both approximations are relatively accurate. In the non-periodic region, the multiple scales approximation outperforms the transseries approximation, while in the transition and 2-periodic region, the transseries approximation is substantially more accurate.

This outcome is sensible; the transseries approximation tracks the contribution of exponentials in the solution, and accurately incorporates them into the solution behaviour. In the non-periodic region, the solution is best represented by an algebraic power series in $\eps$. The multiple scales approach involves calculating this power series to several terms, while our transseries approximation relies only on the leading-order behaviour of this series. In the transition and 2-periodic region, however, these exponential contributions become more significant, and this corresponds to the transseries approximation becoming more accurate than the multiple scales approximation. While the multiple scales approximation is able to capture some of the exponential behaviour, the transseries approximation is able to incorporate several exponential corrections in a straightforward fashion, producing greater accuracy in the solution regions where these corrections play an important role. Furthermore, increasing the accuracy of the transseries approximation in the non-periodic region can be done systematically by including higher corrections in $\varepsilon$.

Finally, we note that there are points in Figure \ref{fig:dyn-error}(b)--(c) where the error appears to drop to zero. This corresponds to a coincidental crossing between the approximation and actual solution occurring at this value of $n$. The crossing may be seen in the example solution from \ref{fig:dyn-approx} at $n \approx 250$. Any solution with a reasonable amount of accuracy will have some value of $n$ where this crossing occurs; this does not provide any added insight into the accuracy of the approximation.

\subsection{4-periodic solution}

\begin{figure}[tb]
\centering
\subfloat[Real part of $B(z)$]{

\includegraphics{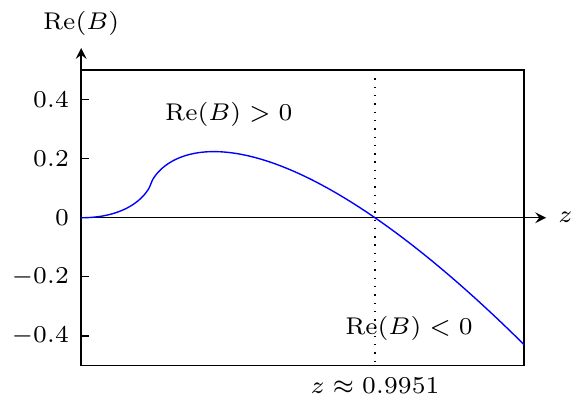}
}
\subfloat[Imaginary part of $B(z)$]{

\includegraphics{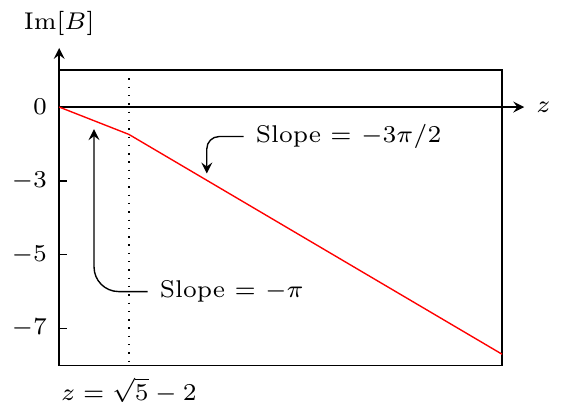}
}
\caption{This figure shows the real and imaginary parts of $B(z)$, corresponding to the exponential weight from \eqref{eq:dyn-4per-weight}. The periodicity of this contribution may be determined by identifying the slope of the imaginary part, corresponding to $\mathrm{Im}[B'(z)]$. For $z < \sqrt{5}-2$, the weight $B(z)$ contains an imaginary term $-\pi \ii$, which corresponds to 2-periodic behaviour. After $z$ exceeds $\sqrt{5}-2$, the slope of the imaginary term changes to $-3\pi\ii/2$, which leads to the appearance of 4-periodic behaviour. This behaviour is not immediately apparent, as the contribution is exponentially small if $\mathrm{Re}[z] > 0$, corresponding to $z < z_0$, where $z_0 \approx 0.9951$. For $z > z_0$, the 4-periodic terms become significant in the solution behaviour. We note that, due to the bifurcation delay, this behaviour is not immediately visibly apparent in the solution; however, a careful analysis of the corresponding transseries terms will identify the transition between 2-periodic and 4-periodic behaviour.}\label{fig:4-per}
\end{figure}

From Figure \ref{fig:dyn-exact}, we see that as $n$ increases, eventually the bifurcation parameter becomes sufficiently large that the solution becomes 4-periodic. As discussed in Section \ref{S:where-4-periodic-behaviour-matters}, this higher periodicity must be encoded in the transseries solution as the weights of new exponential scales. We will not perform a full explicit calculation here, but we will demonstrate that these exponential weights do, in fact, predict the emergence of stable behaviour with higher periodicity. 

We include a new term with transseries parameter $\sigma_1$, and the analysis suggests that it is natural to define a new scaled variable $z = 2 n \eps$. This new contribution to the transseries, denoted $S(z,\eps;\sigma_1)$, is given by
\begin{equation}\label{eq:dyn-4per-transseries}
S(z,\eps;\sigma_1) = \sum_{m=0}^{\infty} \sigma_1^m\e^{-m B(z)/\eps} S_m(z,\eps).
\end{equation}
By adding $S$ as a perturbation to the 2-period solution approximated by \eqref{eq:dynamic-2per-approx} and balancing terms in \eqref{eq:dynamic} in a similar fashion to Section \ref{s:dyn-2per-ansatz}, we obtain an equation for $B(z)$ that gives
\begin{equation}\label{eq:dyn-4per-weight}
B(z) = -\pi \ii  z + z - \frac{z}{2}\log(1 - z(4+z)) + (\sqrt{5}-2)\log\left(\frac{\sqrt{5}-2}{\sqrt{5}-2-z}\right) - (\sqrt{5}+2)\log\left(\frac{\sqrt{5}+2}{\sqrt{5}+2+z}\right),
\end{equation}
where the constant of integration is picked to set $B(0) = 0$ for convenience, though this choice may be absorbed into the parameter $\sigma_1$. The behaviour of $B(z)$ is depicted in Figure \ref{fig:4-per}. There are two significant conclusions that may be drawn from this figure. In Figure \ref{fig:4-per}(b), we see that
\begin{equation}
\mathrm{Im}[B'(z)] = \left\{
        \begin{array}{ll}
            -\pi & \quad z < \sqrt{5}-2 \\
            -3\pi/2 & \quad z > \sqrt{5}-2
        \end{array}
    \right.
    \end{equation}
Noting the format of \eqref{eq:dyn-4per-transseries}, we see that this exponent changes from 2-periodic behaviour to 4-periodic behaviour when crossing the value $z = \sqrt{5}-2$. This change in exponent gives rise to 4 different branches in the solution, and therefore explains the onset of 4-periodic behaviour in the solution to the dynamical logistic equation.

The second important observation is that this 4-periodicity is not immediately apparent in the solution, due to the behaviour of $\mathrm{Re}[z]$. In Figure \ref{fig:4-per}(a), we see that there is a value of $z$, denoted $z_0$ and located at $z_0 \approx 0.9951$, at which the real part of $B(z)$ changes sign from positive to negative. From the form of \eqref{eq:dyn-4per-transseries}, we see that this corresponds to the 4-periodic transseries contribution being exponentially small for $z < z_0$, before growing to have a significant impact on the solution behaviour for $z > z_0$. 

This value of $z_0$ corresponds to 4-periodic behaviour becoming apparent at $n  \approx 0.4975/\eps$. For example, in Figure \ref{fig:dyn-exact}, we would expect that 4-periodic solution to become significant at $n \approx 3455$, which is consistent with the appearance of the second transition region in this image.

A more detailed transseries analysis would permit us to calculate a series approximation for the 4-periodic behaviour; however, as we expected from the transseries approach, a straightforward analysis of the exponential weights in the transseries is sufficient to explain the onset of the higher periodicity, and identify the location in $z$ (and hence, in $n$) where this transition to dominant 4-periodic behaviour takes place.

Finally, we note that the points where the periodicity changes correspond to values of $n$ where the real part of the exponential weights changes sign, or $z_0$ in Figure \ref{fig:4-per}. In asymptotic analysis, this corresponds to the crossing of a curve known as an anti-Stokes curve. This suggests that the Stokes phenomenon plays a role in this system behaviour. In fact, the solution does contain Stokes curves, which are responsible for appearance of exponential factors in the solution; however, finding these Stokes curves requires continuing the solution in the negative-$n$ direction, and was therefore not presented here. Nonetheless, the study the Stokes phenomenon in the dynamic logistic map is an interesting and rich subject which is beyond the scope of the present work.

\section{Discussion}

We have obtained transasymptotic approximations for the solutions to both the standard and slowly varying logistic equation. In each case, we were able not only to reproduce the results calculated in \cite{Hall2016} using multiple scales asymptotic methods, but to go significantly further. 

As, {\it a priori}, transseries methods allow for the straightforward calculation of higher-order exponentials, the transseries approximation was able to represent the solution more accurately than the multiple scales method both during and after the delayed bifurcation, as seen in Figure \ref{fig:dyn-error}; during and after the birfurcation, the initially subdominant exponentials contribute significantly to the solution, so it should be expected that the transseries approximation would be particularly accurate compared to other methods in these regions. 

Furthermore, the transseries approach can still provide a useful approximation when the parameter $\eps$ is not particularly small, as the solution can simply be rescaled to determine the next asymptotic weight. 

We considered the dynamic logistic equation with $\eps > 0$, producing a cascade of delayed bifurcations. If $\eps < 0$, causing the bifurcation parameter to decrease rather than increase, bifurcations appear earlier than the solution stability would suggest, rather than later \cite{Baesens1991a}. A transseries approach could be used in almost identical fashion to the present study in order to approximate these accelerated bifurcations; however, that analysis is beyond the scope of this study.

There are several significant and general advantages to the transseries resummation approach. The first is that the method we have described can be applied in systematic fashion to a wide range of problems, including both discrete and continuous systems. Whilst such advantages have already been seen elsewhere, in the context of the logistic equation it has been instructive to compare this to the multiple scales approach from \cite{Hall2016}, which required the careful comparison of asymptotic terms up to several orders. In order to capture the fast discrete scale, as well as both the inner and outer continuous scales near the bifurcation, asymptotic matching was performed through three scales. 

The transseries approach here was able to reproduce this behaviour by resumming the series in order to ensure that the transasymptotic approximation contained behaviour encoded in the subdominant exponentials; this behaviour contained all of the information found using asymptotic matching methods. Furthermore, improving approximation accuracy by computing more terms of a multiple scales expansion requires comparing the asymptotic behaviour of more terms and checking to determine when the relative dominance of terms changes, while obtaining a more accurate transseries expression simply requires the systematic calculation of further series terms in the transseries. While these calculations can prove challenging, the steps required to obtain the subdominant exponentials, and the associated solution behaviour, follow the same consistent process at each stage which it is applied.

Computing the subdominant exponentials is not valuable simply in that it produces a more accurate approximation. In fact, a second major advantage of the transseries method is that the exponential weights have a significant effect on the system behaviour, and computing just these weights can tell us the form of the solution as parameters in the problem vary. In our analysis of the standard logistic map, we showed that 2-, 4-, and 8-periodic behaviour can be determined simply by carefully studying the subdominant weights. This explains the appearance of higher periodicities in the solution, and suggests that if this process is continued, it can be used to study further bifurcations in the period doubling process. 

In our subsequent analysis of the period doubling cascade found in the slowly varying logistic equation, we were able to predict the onset of 2-periodic and 4-periodic behaviour in the solution, simply by studying the relative size of the exponential weights associated with the 2- and 4-periodic contributions. It would be particularly interesting to continue to investigate how the full period doubling cascade is encoded in the exponential weights of this system, and whether this can provide (at least theoretically) further insight into the period doubling route to chaos.

Finally the full analysis of the movements of exponential contributions between Riemann sheets, seen in the dynamic logistic map, also merits further investigation.  Examples of such phenomena have been observed recently in novel features of aeroacoustic flows \cite{stone2017,stone2018}. Initial explorations appear to suggest this is commonly found in other physical and mathematical contexts.

\paragraph*{Acknowledgements} IA has been supported by the UK EPSRC Early Career Fellowship EP/S004076/1, and the FCT-Portugal grant PTDC/MAT-OUT/28784/2017. CJL has been supported by Australian Resaerch Council Discovery Project \#190101190. DH has been supported by
the presidential scholarship of the University of Southampton. The authors thank Dr Cameron Hall for helpful discussions.

\bibliographystyle{plain}
\bibliography{doublerefs.bib}

\appendix

\section{Explicit transseries terms}\label{A:logistic-transseries}

In \eqref{eq:logistic-R-newdef}, the transseries for $R(\tau_0,\eps)$ is written in terms of a base approximation $R_0(\eps)$, a sum of odd terms in $\tau_0$, denoted $\Omega_{o,k}$, and a sum of even terms in $\tau_0$, denoted $\Omega_{e,k}$. We further simplify this by writing $R(\tau_0,\eps;\sigma_0) = R(\eps) + S(\tau_0,\eps)$, where
\begin{equation}
S(\tau_0,\eps) =  \sqrt{\eps} \sum_{k=0}^{\infty}\eps^{k} \Omega_{o,k}(\tau_0) + \eps \sum_{k=0}^{\infty} \eps^k \Omega_{e,k}(\tau_0).\label{eq:logistic-delR}
\end{equation}
We note that $\tau_0(x+\eps) = -(1+\eps)\tau_0(x)$. Consequently,
\begin{equation}
S(\tau_0(x+\eps),\eps) = S(-(1+\eps)\tau_0(x),\eps).\label{eq:tau_x_plus_eps}
\end{equation}
Applying \eqref{eq:logistic-delR} and \eqref{eq:tau_x_plus_eps} to the logistic equation \eqref{eq:logistic-finite-diff-constant} gives
\begin{equation}
S(-(1+\varepsilon)\tau_0,\eps)=-(1+\varepsilon)S(\tau_0,\eps)-(3+\varepsilon)S(\tau_0,\eps)^{2}.\label{eq:logistic-R-eqn}
\end{equation}
Expanding the left-hand side of this expression as a Taylor series in $\varepsilon$ gives
\begin{align}
\nonumber S(-(1+\eps)\tau_0,\eps) &= \sum_{j=0}^{\infty} \frac{(-\tau_0 \eps)^j}{j!}R^{(j)}(-\tau_0) \\
&= -\sqrt{\eps}\sum_{m=0}^{\infty}\eps^m \sum_{k=0}^m \frac{\tau_0^k}{k!}\Omega^{(k)}_{o,m-k}(\tau_0) + \eps \sum_{m=1}^{\infty}\eps^{n}\sum_{k=0}^{m-1}\frac{\tau_0^k}{k!}\Omega^{(k)}_{e,n-1-k}(\tau_0),
\end{align}
where we used the fact that $\Omega_{o,k}$ and $\Omega_{e,k}$ are odd and even in $\tau_0$ respectively. The remaining expansions in \eqref{eq:logistic-R-eqn} may be obtained by substitution of \eqref{eq:logistic-delR} into \eqref{eq:logistic-R-eqn}. It is straightforward to show that
\begin{align}
R(\tau_0,\eps)^2 = 3\sum_{m=1}^{\infty}\varepsilon^{m} & \sum_{k=0}^{m-1}\Omega_{o,k}(\tau_0)\Omega_{o,m-1-k}(\tau_0)\\
& +\sqrt{\varepsilon}\sum_{m=1}^{\infty}\varepsilon^{m}\,\sum_{k=0}^{m-1}\Omega_{o,k}(\tau_0)\Omega_{e,m-1-k}(\tau_0)
 +\sum_{m=2}^{\infty}\varepsilon^{m}\,\sum_{k=0}^{m-2}\Omega_{e,k}(\tau_0)\Omega_{e,m-2-k}(\tau_0).
\end{align}
These expansions may now be used to equate powers of $\eps$ and obtain the expressions given in \eqref{eq:static-omega0} and \eqref{eq:static-omega1}.

\section{Initial Condition for 4-Periodic Equation}\label{A.2}

In order to calculate the initial condition for the 4-periodic problem, we first recall that $\hat{R}(x,\eps)$ was derived in order to satisfy the initial condition for small $\epsilon$. The 4-periodic solution arises for $\eps > -2 + \sqrt{6}$, or $\eta > 0$. Hence, we determine the initial condition by perturbing around the leading-order behaviour of $\hat{R}(x,\eps)$, which is initially 2-periodic for the parameter regime under consideration. We then determine $\sigma_1$ by matching with the initial condition in the limit that $\eta \rightarrow 0$.

We first obtain stable 2-periodic behaviour of $R(x,\eps)$ from \eqref{eq:2per-manifold}, letting $x = 0$ in order to describe the initial state. This expression may be written in terms of $\eta$, to allow a small $\eta$ expansion in this limit. This gives
\begin{align}\nonumber
\hat{R}(0,\eps) = \frac{4 + \eps + \sqrt{\eps(4+\eps)}}{2(3+\eps)} \sim \frac{1}{5}\left(2 - \sqrt{3} + \sqrt{2 + \sqrt{3}}\right) + \frac{\eta}{50}\left(3\sqrt{2} - 16\sqrt{3} - 7 \sqrt{6} + 12\right)& \\
+ \frac{3\eta^2}{250}\left(-47 \sqrt{2} + 84\sqrt{3} + 18\sqrt{6} - 38\right)& + \mathcal{O}(\eta^3).
\end{align}
Setting $x = 0$, letting $\sigma_1 = \sigma_{1,0} + \eta\sigma_{1,1} + \ldots$, and expanding $S$ in powers of $\eta$ gives
\begin{align}
S(0,\eps) \sim \eta \sigma_{1,0} + \eta^2\left(-\tfrac{5}{12}(14-7\sqrt{2}-4\sqrt{3}+4\sqrt{6})\sigma_{1,0}^3 + (3\sqrt{3}-\sqrt{6}+\tfrac{3}{2}\sqrt{2}-1)\sigma_{1,1}^2 + \sigma_1\right) + \mathcal{O}(\eta^3)
\end{align}
To determine the appropriate initial condition, we fix the case for $\eta = 0$, which gives
\begin{equation}
R(0,\eps) =  \frac{1}{5}\left(2 - \sqrt{3} + \sqrt{2 + \sqrt{3}}\right).
\end{equation}
By setting $R(0,\eps) = \hat{R}(0,\eps) + S(0,\eps)$, and matching powers of $\eta$, we can obtain
\begin{equation}
\sigma_1 = - \frac{1}{50}\left(3\sqrt{2} - 16\sqrt{3} - 7 \sqrt{6} + 12\right) + \frac{\eta}{500}\left(297\sqrt{2} - 709\sqrt{3} - 189\sqrt{6} + 399\right) + \mathcal{O}(\eta^2),
\end{equation}
This is sufficient information to approximate the solution using the transseries behaviour, although it is straightforward to continue this process to obtain higher corrections for $\sigma_1$.

\end{document}